\documentclass[10pt,letterpaper]{amsart}%
\usepackage[utf8]{inputenc}%
\usepackage{framed,multirow}%
\usepackage{amssymb,amsmath,amsthm,amstext,amscd}%
\usepackage{graphicx}%
\usepackage{subcaption}%
\usepackage[dvipsnames]{xcolor}%
\usepackage[per-mode=symbol]{siunitx}%
\DeclareSIUnit\Length{L}%
\DeclareSIUnit\Time{T}%
\DeclareSIUnit\Mass{M}%
\DeclareSIUnit\liter{\ell}%
\DeclareSIUnit\nats{nats}%
\usepackage{multirow,bigdelim,booktabs,colortbl}%
\usepackage{enumitem}%
\usepackage[ruled, vlined]{algorithm2e}%
\usepackage{lineno}

\usepackage{todonotes}%

\usepackage{hyperref}%
\usepackage{cleveref}%
\usepackage[normalem]{ulem}
\theoremstyle{plain}%
\newtheoremstyle{example}%
{}%
{}%
{}% Body font
{}% Indent amount (empty = no indent, \parindent = para indent)
{\it}% Thm head font
{}% Punctuation after thm head
{\newline}% Space after thm head (\newline = linebreak)
{\thmname{#1}\thmnumber{ #2}: \thmnote{#3}.}% Thm head spec
\newtheoremstyle{remark}%
{}%
{}%
{}% Body font
{}% Indent amount (empty = no indent, \parindent = para indent)
{\it}% Thm head font
{:}% Punctuation after thm head
{ }% Space after thm head (\newline = linebreak)
{\thmname{#1}\thmnumber{ #2}}% Thm head spec
\theoremstyle{remark}%
\newtheorem{remark}{Remark}%
\crefname{equation}{}{}%
\crefname{remark}{Remark}{Remarks}%
\crefname{theorem}{Theorem}{Theorems}%
\crefname{lemma}{Lemma}{Lemmas}%
\crefname{corollary}{Corollary}{Corollaries}%
\crefname{figure}{Figure}{Figures}%
\crefname{algorithm}{Algorithm}{Algorithms}%
\crefname{section}{Section}{Sections}%
\crefname{subsection}{Section}{Sections}%
\crefname{table}{Table}{Tables}%

\DeclareMathOperator{\E}{\mathbf{E}}%
\newcommand{\Pb}{P}%
\newcommand{\dd}{\mathrm{d}}%
\DeclareMathOperator{\RE}{\mathcal{R}}%
\newcommand{\unif}{\mathrm{unif}}%
\newcommand{\given}{\mid}%
\newcommand{\midbars}{\;\|\;}%
\newcommand{\ind}{\mathrel{\perp}}%
\newcommand{\defeq}{\mathrel{\mathop:}=}%
\newcommand{\eqdef}{=\mathrel{\mathop:}}%
\newcommand{\Deff}{\mathbf{D}_{\mathrm{eff}}}%
\newcommand{\DL}{{D_{L}}}%
\newcommand{\DT}{{D_{T}}}%
\newcommand{\geff}{\gamma_{\mathrm{eff}}}%
\newcommand{\B}[1]{\boldsymbol{#1}} \newcommand{\vol}{\mathcal{V}}%
\newcommand{\por}{\mathcal{P}}%
\newcommand{\sol}{\mathcal{S}}%
\newcommand{\infc}{\Gamma}%
\newcommand{\uvol}{{\breve{\mathcal{V}}}}%
\newcommand{\usol}{{\breve{\mathcal{S}}}}%
\newcommand{\upor}{{\breve{\mathcal{P}}}}%
\newcommand{\uinfc}{{\breve{\Gamma}}}%
\newcommand{\iid}{{\scshape iid}}%
\usepackage[margin=1in]{geometry}%

\title[Bayesian network PDEs for multiscale porous media]{Causality and Bayesian network PDEs for multiscale
  representations of porous~media}%

\author[K.\ Um]{Kimoon Um$^*$}%
\address[K.\ Um]{Department of Energy Resources Engineering, Stanford
  University, 367 Panama Street, Stanford, CA 94305, USA.}
\email{kimoon.um@gmail.com}%

\author[E.\,J.\ Hall]{Eric J.\ Hall$^*$}%
\address[E.\,J.\ Hall]{Chair of Mathematics for Uncertainty Quantification, RWTH Aachen University, Kackertstraße 9 C, Aachen 52072, Germany.}
\email{hall@uq.rwth-aachen.de}%
\thanks{$^*$Both authors contributed equally to this work.}
  
\author[M.\,A.\ Katsoulakis]{Markos A.\ Katsoulakis}%
\address[M.\,A.\ Katsoulakis]{Department of Mathematics and Statistics, University of
  Massachusetts Amherst, 710 N. Pleasant Street, Amherst, MA 01003,
  USA.}
\email{markos@math.umass.edu}%
  
\author[D.\,M.\ Tartakovsky]{Daniel M.\ Tartakovsky}%
\address[D.\,M.\ Tartakovsky]{Department of Energy Resources Engineering, Stanford
  University, 367 Panama Street, Stanford, CA 94305, USA.}
\email[Corresponding author]{tartakovsky@stanford.edu}%

\crefname{equation}{}{}%
\crefname{remark}{Remark}{Remarks}%
\crefname{theorem}{Theorem}{Theorems}%
\crefname{lemma}{Lemma}{Lemmas}%
\crefname{corollary}{Corollary}{Corollaries}%
\crefname{figure}{Figure}{Figures}%
  
\begin{document}%
\maketitle%

\begin{abstract}
  Microscopic (pore-scale) properties of porous media affect and often
  determine their macroscopic (continuum- or Darcy-scale)
  counterparts. Understanding the relationship between processes on
  these two scales is essential to both the derivation of macroscopic
  models of, e.g., transport phenomena in natural porous media, and
  the design of novel materials, e.g., for energy storage. Most
  microscopic properties exhibit complex statistical correlations and
  geometric constraints, which presents challenges for the estimation
  of macroscopic quantities of interest (QoIs), e.g., in the context
  of global sensitivity analysis (GSA) of macroscopic QoIs with
  respect to microscopic material properties. We present a systematic
  way of building correlations into stochastic multiscale models
  through Bayesian networks. This allows us to construct the joint
  probability density function (PDF) of model parameters through
  causal relationships that emulate engineering processes, e.g., the
  design of hierarchical nanoporous materials. Such PDFs also serve as
  input for the forward propagation of parametric uncertainty; our
  findings indicate that the inclusion of causal relationships impacts
  predictions of macroscopic QoIs. To assess the impact of
  correlations and causal relationships between microscopic parameters
  on macroscopic material properties, we use a moment-independent GSA
  based on the differential mutual information. Our GSA accounts for
  the correlated inputs and complex non-Gaussian QoIs. The global
  sensitivity indices are used to rank the effect of uncertainty in
  microscopic parameters on macroscopic QoIs, to quantify the impact
  of causality on the multiscale model's predictions, and to provide
  physical interpretations of these results for hierarchical
  nanoporous materials.
\end{abstract}


\section{Introduction}%
\label{sec:introduction}%

Understanding statistical and causal relations between
properties/model parameters at various scales is essential for
science-based predictions in general and for forecasts of transport
phenomena in porous media in particular. For example, the design of
materials for energy storage devices aims to optimize macroscopic
material properties (quantities of interest or QoIs), such as
effective diffusion coefficient and capacitance, through engineered
pore structures~\cite{ZhangEtAl:2015od,ZhangTartakovsky:2017od}.
Quantification of both uncertainty in predictions of these macroscopic
QoIs and their sensitivity to variability and uncertainty in
microscopic features are crucial for informing such decision tasks as
optimal experimental design and reliability engineering.

The simulation-assisted approach to the optimal design of porous
meta-materials takes advantage of the availability of microscopic
(pore-scale) and macroscopic (continuum or Darcy-scale) models, as
well as of a bridge between the two provided by various upscaling
techniques~\cite{ZhangEtAl:2015od,ZhangTartakovsky:2017od,Ling:2016di}.
Such a bridge also facilitates analysis of both uncertainty
propagation from the microscopic scale to the macroscopic scale and
sensitivity of macroscopic material properties to the microscopic
ones~\cite{UmZhangKatsoulakisEtAl:2017aa}. Predictions resulting from
this multiscale approach can be made more robust by incorporating
information about correlations and causal relationships between scales
and within a single scale. Causality stems for example from physical,
chemical, and/or engineering design constraints. Our primary objective
is to bring a Bayesian network perspective to the incorporation of
causality into modeling process.

Bayesian networks are a special class of hierarchical probabilistic
graphical models (PGMs) with a directed acyclic graph structure that
represent causal relationships among random
variables~\cite{pearl1988probabilistic,pearl2014probabilistic,
  KollerFriedman:2009gm}. %Wasserman:2013as,Bishop:2006ml}.
Bayesian networks and, more generally, PGMs provide a rich framework
for encoding distributions over large, complex domains of interacting
random variables that can include causal relationships and expert
knowledge.

In our application, they provide a coherent framework for representing
causal relationships both between problem scales and among the space
of parameters representing pore-scale features. Incorporating Bayesian
networks into random PDEs is a novel approach to the modeling of
multiscale porous media; it breaks down the stochastic modeling and
statistical inference tasks into smaller, controllable parts enabling
us to
\begin{enumerate}[label=(\roman*)]
\item build systematically informed parameter priors that include
  physical constraints and/or correlations;
\item mirror engineering processes related to the design of
  hierarchical nanoporous media networks;
\item construct Bayesian network (random) Darcy-scale PDE models
  informed by (possibly uncertain) pore-scale data, parameters, and
  constraints; and
\item carry out global sensitivity analysis (GSA) and uncertainty
  quantification (UQ).
\end{enumerate}

This framework for incorporating causal relationships through
structured priors demands GSA tools that differ from the standard
variance-based GSA methods. These typically assume unstructured (i.e.,
mutually independent) priors and are neither easy to interpret nor
cheap to compute for correlated
inputs~\cite{MaraEtAl:2015np,IoossPrieur:2018se}. In our application,
causal relationships exist not just between parameters but also
between scales. The latter is important since our predictive PDFs are
not necessarily Gaussian and/or do not have a known analytical form.
This suggests moving away from a fixed number of moments to a
moment-independent quantity such as mutual information. We employ a
moment-independent GSA relying on mutual
information~\cite{CoverThomas:2006in,Soofi:1994if} and empirical
distributions acquired through simulations. We demonstrate that
differential mutual information provides a measure of input effects
that is suited to tackling the twin challenges of structured or
correlated inputs and non-Gaussian QoIs.

Design of (nano)porous metamaterials for energy storage provides an
ideal setting to illustrate the Bayesian network PDE approach.
Macroscopic material properties are dependent on a set of microscopic
parameters characterizing the pore geometry, e.g., pore radius or pore
connectivity. These microscopic parameters are typically correlated
due to the presence of geometric and topological constraints and
uncertain due to natural and/or manufacturing variability. This
setting gives rise to a number of theoretical and practical questions:
How does uncertainty in microscopic properties (quantified, e.g., in
terms of a pore-size distribution) propagate to the macroscopic scale
(expressed in terms of the PDF of, e.g., the effective diffusion
coefficient)? How sensitive are a material's macroscopic properties to
its microscopic counterparts? Etc. Our computational strategy, which
makes exhaustive sampling for prediction and uncertainty
quantification feasible, has three ingredients: i) Rosenblatt
transforms to decorrelate inputs for non-intrusive scientific
computing, ii) generalized polynomial chaos expansions obtained using
the \texttt{DAKOTA} software~\cite{DAKOTA:2009man}, and iii) kernel
density estimation techniques.

In \cref{sec:pore-to-darcy-model}, we formulate a macroscopic
(Darcy-scale) model of reactive transport in hierarchical nanoporous
media; the model parameters are expressed in terms of microscopic
(pore-scale) material properties by means of
homogenization~\cite{UmZhangKatsoulakisEtAl:2017aa}, which facilitates
multiscale UQ and GSA. \cref{sec:prob-graph-model} contains a
description of our Bayesian network-based approach to linking the
components of this model across the two scales.
\cref{sec:quant-uncert-darcy} contains details of its implementation
highlighting the use of the inverse Rosenblatt transform to
non-intrusive utilization of existing software, such as
\texttt{DAKOTA}. This section also collates results of our numerical
experiments, which demonstrate the importance of causality. In
\cref{sec:global-sa}, we adopt sensitivity indices based on
differential mutual information and provide ranking for the impact of
uncertainty in the microscopic parameters on uncertainty in their
macroscopic counterparts. Our examples illustrate how the inclusion of
causal relationships encoding structural constraints provides rankings
more consistent with the physics anticipated for a simple hierarchical
nanoporous material. In \cref{sec:alternative-prob-models}, we present
an alternative Bayesian network to highlight the method's flexibility
and ability to mirror distinct engineering and design processes. Major
conclusions drawn from our study are summarized in
\cref{sec:conclusions}.

\section{Models of flow and transport in nanoporous materials}%
\label{sec:pore-to-darcy-model}%

\begin{figure}[!h] \centering
  \includegraphics[width=.6\textwidth]{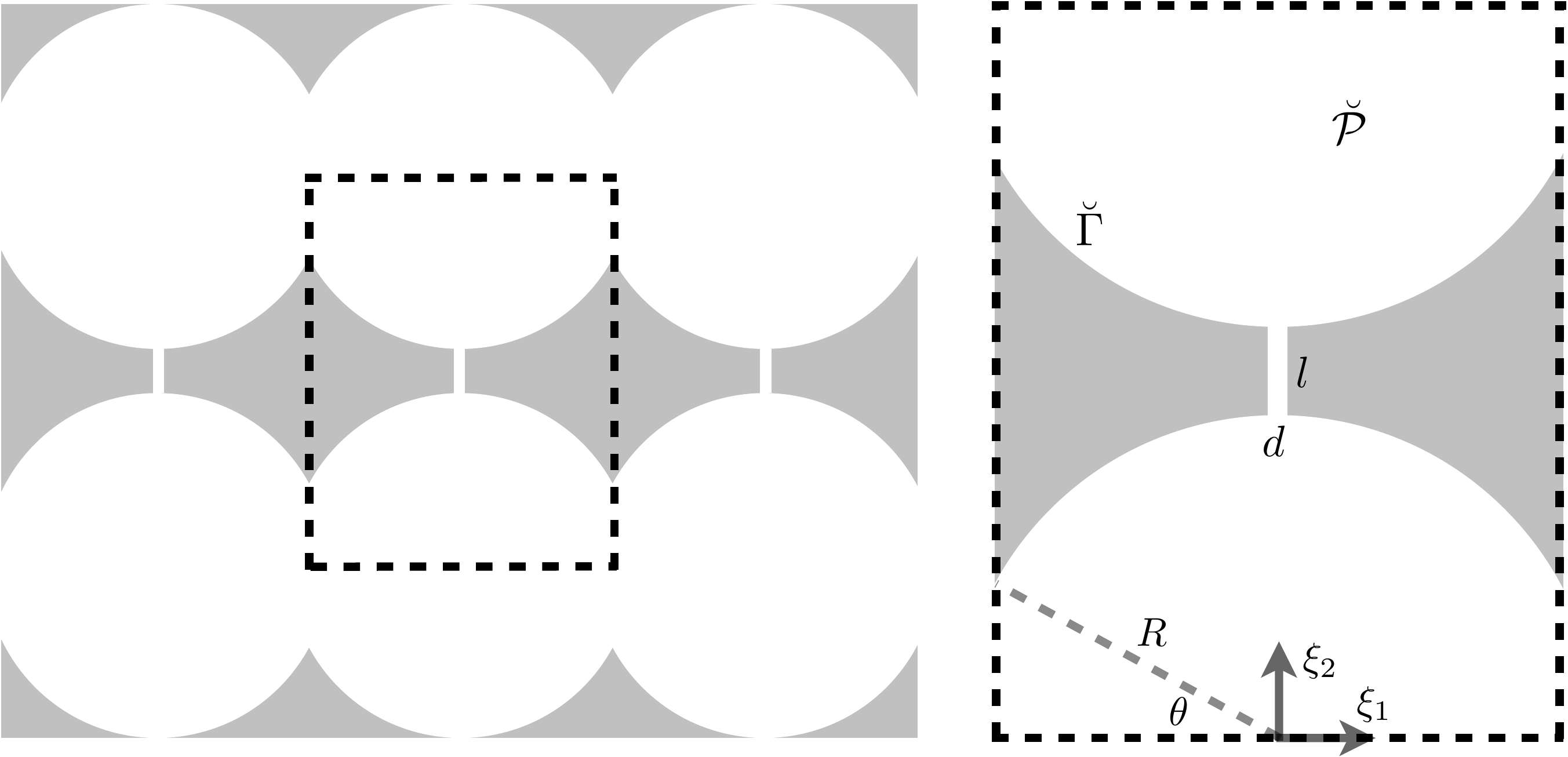}
  \caption{A hierarchical nanoporous
    material~\cite{UmZhangKatsoulakisEtAl:2017aa} exhibiting
    horizontally oriented nanotunnels through mesopores connected by a
    series of vertically oriented nanotubes. The porous media volume
    $\vol$ (left) consists of a periodic arrangement of unit cells
    $\uvol$ (right) with pore space $\upor$ and fluid-solid interface
    $\uinfc$. The parameters %$\B{p} = $
    $\{R, \theta, d, l\}$ describing the nanopore features are
    constrained by the geometry of $\uvol$. }
  \label{fig:pore-structure}
\end{figure}

A volume $\vol = \por \cup \sol$ of a hierarchical porous material is
comprised of a fluid-filled pore space $\por$ and solid matrix $\sol$,
with a (multi-connected) fluid-solid interface denoted by
$\infc = \por \cap \sol$. To mimic a manufacturing process and to
make subsequent use of the homogenization theory, we assume that the
volume $\vol$ consists of a periodic arrangement of unit cells
$\uvol = \upor \cap \usol$ with pore space $\upor \subset \por$, solid
matrix $\usol \subset \sol$, and fluid-solid interface
$\uinfc = \upor \cap \usol$. For example, the hierarchical nanoporous
material in \cref{fig:pore-structure} consists of mesopores that are
connected longitudinally (horizontally) through nanotunnels and
transversely (vertically) by a series of nanotubes. These features are
described by a set of parameters $\{R, \theta, d, l\}$, where $R$ is
the mesopore radius; $\theta$ is the angle of overlap between adjacent
mesopores in a nanotunnel; and $d$ and $l$ are the diameter and length
of the nanotubes, which serve as nano-bridges between adjacent
mesopores/nanotunnels.

The design of novel materials calls for a systematic analysis of the
sensitivity of desired macroscopic (Darcy-scale) properties to
imperfections (natural variability) in microscopic (pore-scale)
parameters and/or their distributions.
Following~\cite{UmZhangKatsoulakisEtAl:2017aa}, we use the
homogenization theory to map uncertainty in the microscopic parameters
and processes to their macroscopic counterparts. Sources of
uncertainty, as well as representations of randomness in the
microscopic and macroscopic models are described below.

\subsection{Pore-scale model}%
\label{sec:pore-scale-model}%

At the pore-scale, the evolution of a solute concentration
$c(\B{x}, t)$ (\si{\mol\per\liter\cubed}), at point $\B{x} \in \por$
and time $t > 0$, is governed by the evolution equation
\begin{equation}
  \label{eq:pore-solute-evolution}
  \frac{\partial c}{\partial t} = \nabla \cdot (D \nabla c)\,, 
  \qquad \B{x} \in \por\,, \qquad t > 0\,, 
\end{equation}
where $D(\B{x})$ (\si{\Length\squared\per\Time}), $\B{x} \in \por$, is
the pore-scale diffusion coefficient. The spatial variability of $D$
allows for Fickian diffusion through mesopores and Knudson diffusion
through nanotubes. This equation is subject to the uniform initial
condition
\begin{equation*}
  \label{eq:pore-solute-init-cond}
  c(\B{x},0) = c_{\mathrm{in}}, \qquad \B{x} \in \por\,,
\end{equation*}
and the boundary condition
\begin{equation*}
  \label{eq:pore-solute-bdry-cond}
  -D\B{n} \cdot\nabla c = q_{\mathrm{m}}\frac{\partial s}{\partial t} \,, 
  \qquad \B{x} \in \infc \,, \qquad t > 0 \,,
\end{equation*}
where $q_{\mathrm{m}}$ and $s(\B{x},t)$ are related to the sorption
properties of the material surface $\infc$. Specifically,
$q(\B{x}, t) = q_{\mathrm{m}} \cdot s(\B{x},t)$ is the adsorption
amount per unit area of $\infc$ (\si{\mol\per\square\liter}),
$q_{\mathrm{m}}$ (\si{\mol\per\liter\squared}) is the maximal
adsorption amount, and $s(\B{x},t)$ is the fractional coverage of
$\infc$. The fractional coverage is assumed to follow Lagergren's
pseudo-first-order rate equation,
\begin{equation*}
  \frac{\dd s}{\dd t} = \gamma(s_{\mathrm{eq}} -s )\,,
\end{equation*}
where $\gamma$ (\si{\per\Time}) is the adsorption rate constant and
the equilibrium adsorption coverage fraction $s_{\mathrm{eq}}$
satisfied Langmuir's adsorption isotherm,
\begin{equation}
  \label{eq:Langmuir-adsoprtion-isotherm}
  s_{\mathrm{eq}} = \frac{Kc}{1+Kc}\,,
\end{equation}
with the adsorption equilibrium constant $K$
(\si{\liter\cubed\per\mol}).

\subsection{Darcy-scale model}%
\label{sec:darcy-scale-model}%

At the macroscopic scale, the volume-averaged solute concentration
$u$,
\begin{equation*}
  \label{eq:Darcy-concentration}
  u(\B{x},t) \defeq 
  \frac{1}{\|\uvol\|} \int_{\uvol(\B{x})} c(\B{\xi}, t) \dd \B{\xi} 
  = \frac{1}{\|\uvol\|} \int_{\upor(\B{x})} c(\B{\xi}, t) \dd \B{\xi} 
  = \frac{\phi}{\|\upor\|} \int_{\upor(\B{x})} c(\B{\xi}, t) \dd \B{\xi}\,, 
  \qquad \B{x} \in \vol\,,
\end{equation*}
treats a porous material as a continuum. Here, $\|\cdot\|$ denotes the
volume of a domain and
$\phi \defeq \|\por\|/\|\vol\| = \|\upor\|/\|\uvol\|$ is the material
porosity. Using homogenization via multiple-scale
expansions~\cite{ZhangEtAl:2015od}, one can show that $u$ satisfies a
reaction-diffusion equation
\begin{equation}
  \label{eq:Darcy-solute-evolution}
  \phi \frac{\partial u}{\partial t} = 
  \nabla\cdot(\Deff\nabla u) - \phi q_{\mathrm{m}}\geff\frac{Ku}{1+Ku}, \qquad \B{x} \in \vol.
\end{equation}
The effective diffusion coefficient $\Deff$ and the effective rate
constant $\geff$ are random, stemming from uncertainty in pore-scale
structures and processes that is propagated by the homogenization map.
Specifically, the effective rate constant $\geff$ (\si{\per\Length})
is computed as
\begin{equation}
  \label{eq:geff}
  \geff = \frac{\|\uinfc\|}{\|\upor\|}\,,
\end{equation}
i.e., is defined solely by the pore geometry; and the effective
diffusion coefficient $\Deff$ (L$^2$/T), a second rank tensor, depends
on both the pore geometry and the pore-scale processes. It is computed
in terms of a closure variable $\B{\chi}$ as
\begin{equation}
  \label{eq:Deff}
  \Deff = \frac{1}{\|\uvol\|} \int_{\upor} (\B{I} + \nabla_{\B{\xi}}\B{\chi}^\top)\dd\B{\xi}\,,
\end{equation}
where $\B{I}$ is the identity matrix. The closure variable
$\B{\chi}(\B{\xi})$ is a $\uvol$-periodic vector defined on $\upor$,
which satisfies the Laplace equation
\begin{equation}
  \label{eq:closure-equation}
  \nabla_{\B{\xi}} \cdot(D\nabla_{\B{\xi}}\B{\chi}) = \B{0}\,, 
  \qquad \B{\xi} \in \upor\,,
\end{equation}
subject to the normalizing condition
\begin{equation}
  \label{eq:closure-norm-cond}
  \langle \B{\chi} \rangle 
  \defeq \frac{1}{\|\uvol\|} \int_\upor\B{\chi}(\B{\xi}) \dd\B{\xi} 
  = \B{0} \,,
\end{equation}
the boundary condition along the fluid-solid segments $\uinfc$,
\begin{equation}
  \label{eq:closure-bdry-cond}
  \B{n} \cdot \nabla_{\B{\xi}}\B{\chi} = - \B{n}\cdot \B{I}\,, 
  \qquad \B{\xi} \in \uinfc\,,
\end{equation}
and $\uvol$-periodic boundary condition on the remaining fluid
segments of the boundary of $\upor$. For the hierarchical nanoporous
material in \cref{fig:pore-structure}, these auxiliary conditions
reduce to
\begin{equation}
  \label{eq:eg-closure-cond-chi1}
  \chi_1(-a,\xi_2) = \chi_1(a,\xi_2) = 0\,, 
  \qquad \frac{\partial\chi_1}{\partial\xi_2}(\xi_1,0) 
  = \frac{\partial\chi_1}{\partial\xi_2} (\xi_1,b) = 0\,,
\end{equation}
and
\begin{equation}
  \label{eq:eg-closure-cond-chi2}
  \chi_2(\xi_1, 0) = \chi_2(\xi_1, b) = 0\,, 
  \qquad \frac{\partial\chi_2}{\partial\xi_1}(-a, \xi_2) 
  = \frac{\partial\chi_2}{\partial\xi_1} (a, \xi_2) = 0\,,
\end{equation}
where
\[a = R\cos\theta \quad\text{and}\quad b = 2R\cos(\sin^{-1}
  (\frac{d}{2R})) + l = 2R\sqrt{1-\frac{d^2}{4R^2}} + l
\]
are, respectively, the width and height of the unit cell $\uvol$. In
the next section, we introduce a representation for uncertainty in the
pore- and Darcy-scale models that will enable us to compartmentalize
various stochastic modeling and statistical inference tasks.

\section{Bayesian network formulation for random PDE models of
  multiscale porous media}
\label{sec:prob-graph-model}

\begin{figure}[!h]
  \centering \includegraphics[width=.4\textwidth]{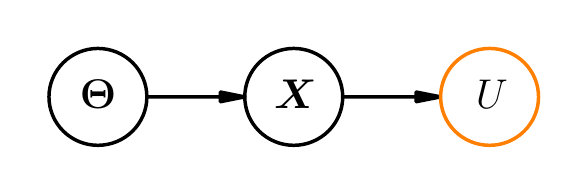}
  \caption{A Bayesian network describing the components of the full
    statistical model $P$, in \cref{eq:statistical-model}, for the
    multiscale porous media system takes into account the joint PDF
    $\Pb(\B{\Theta})$ of the input parameters, the PDF
    $\Pb(\B{X} \given \B{\Theta})$ of the upscaling variable that maps
    pore-scale properties to Darcy-scale variables, and the PDF
    $\Pb(U \given \B{X})$ related to Darcy-scale QoIs. This figure and
    other Bayesian networks are produced using \cite{scipy:daft}.}
  \label{fig:general-model}
\end{figure}

A Bayesian network is a directed graph structure, in which each node
represents a random variable with an associated PDF and each edge
encodes a causal relationship \cite{KollerFriedman:2009gm}. A Bayesian
network PDE incorporates these structured probabilistic models into a
forward physical model, and allows one to capture causal relationships
in a systematic way. The key components of the full statistical model
include:
\begin{enumerate}[label=(\roman*)]
\item inputs $\B{\Theta} = \{R, \theta, d, l\}$, a random vector
  related to the parameters describing pore-scale features in
  \cref{fig:pore-structure};
\item upscaling variable $\B{X}$, a random vector related to the
  closure
  equations~\eqref{eq:closure-equation}--\eqref{eq:closure-bdry-cond};
  and
\item QoI $U$, a random macroscopic quantity, such as macroscopic
  parameters $\Deff$ and $\geff$ in \cref{eq:geff,eq:Deff} or a
  functional of $u$ in \cref{eq:Darcy-solute-evolution}.
\end{enumerate}
Causal relationships arise naturally among the problem scales and
hence these components: the PDF of $U$ depends on the PDF of $\B{X}$
since, for example, $u$ and $\Deff$ depend on the closure variable
$\B{\chi}$. In turn, the PDF of $\B{X}$ depends on the PDF of
$\B{\Theta}$ since $\B{\chi}$ depends on the pore-scale parameter
values, e.g.,
via~\cref{eq:eg-closure-cond-chi1,eq:eg-closure-cond-chi2}. The
Bayesian network in \cref{fig:general-model} describes these causal
relationships for the full statistical model. The directed network
structure indicates that there is only one-way communication between
the different scales.

More precisely, the Bayesian network in \cref{fig:general-model}
encodes conditional relationships among the PDFs for $\B{\Theta}$,
$\B{X}$, and $U$. We denote the joint PDF of the input parameters by
$\Pb(\B{\Theta})$. For simplicity we assume the statistical models for
$\B{X}$ and $U$ to be known, that is, uncertainty enters only through
the parameters since we have fixed both the form of equations for
$\B{\chi}$ and $u$ in
\eqref{eq:Darcy-solute-evolution}--\eqref{eq:closure-bdry-cond} and
the numerical methods for approximating them. Under this assumption
the conditional PDF of $\B{X}$ given a sample $\B{\Theta}$ %$ = \B{p}$
is the Dirac delta function
\begin{equation}
  \label{eq:dist-X}
  \Pb(\B{X} \given \B{\Theta}) 
  = \delta_{\B{X}} (\B{X} - \B{\chi}(\B{\xi}; \B{\Theta})).
\end{equation}
Similarly, the conditional PDF of $U$ given a sample
$\B{X} = \B{\chi}$ is
\begin{equation}
  \label{eq:dist-U}
  \Pb(U \given \B{X}) 
  = \delta_{U} (U - u(\B{x},t; \B{X})).
\end{equation}
Then the full statistical model $P$ is given by
\begin{equation}
  \label{eq:statistical-model}
  \begin{split}
    P \defeq \Pb(\B{\Theta}, \B{X}, U) = \Pb(U \given \B{X}) \;
    \Pb(\B{X} \given \B{\Theta})\; \Pb(\B{\Theta}) = \delta_{U} (U -
    u(\B{x},t; \B{X})) \; \delta_{\B{X}} (\B{X} - \B{\chi}(\B{\xi},t;
    \B{\Theta}))\; \Pb(\B{\Theta})\,.
  \end{split}
\end{equation}
In the remainder of this section, we discuss expanded Bayesian
networks for representing causal relationships between parameters that
allow us to encode correlations among pore scale features. In
\cref{sec:uniform-priors}, we consider the assumption of independent
priors for the pore-scale parameters and contrast this in section
\cref{sec:correlations-from-constraints} with causality arising from
natural structural constraints encountered in engineering design.

\begin{remark}
  Widely adopted approaches to uncertainty quantification, e.g., Monte
  Carlo methods, involve strategies for generating surrogates of the
  model $P$ in \cref{eq:statistical-model}. The direct application of
  these traditional UQ methods require the form of $P$ to be known. In
  the sequel we introduce an approach relying on Rosenblatt
  transformations, generalized polynomial chaos expansions, and the
  popular UQ software package \texttt{DAKOTA} that makes sampling $P$
  feasible.
\end{remark}

\subsection{Independent uniform priors}%
\label{sec:uniform-priors}%

\begin{figure}[htbp]
  \centering \includegraphics[width=.4\textwidth]{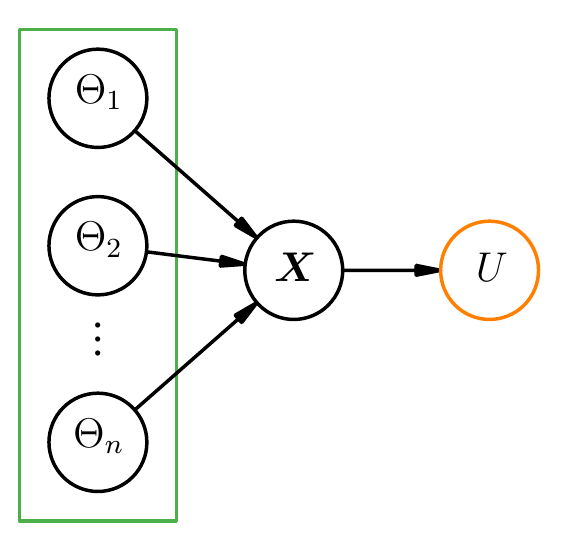}
  \caption{A Bayesian network describing the components of the full
    statistical model under the assumption of independent priors on
    pore-scale features $\B{\Theta} = (\Theta_1, \dots, \Theta_n)$.
    The flat structure of the $\B{\Theta}$ component in the model
    above contrasts with the rich structure of the Bayesian network in
    \cref{fig:conditional-d-and-l} that captures causal relationships
    among the pore-scale features in order to ensure sampling
    geometries consistent with the hierarchical nanoporous material in
    \cref{fig:pore-structure} over the physically relevant
    hyperparameter ranges in \cref{tab:extended-range}.}
  \label{fig:general-model-indpriors}
\end{figure}

As a first simple case we revisit the
analysis~\cite{UmZhangKatsoulakisEtAl:2017aa} of the hierarchical
nanopore geometry in \cref{fig:pore-structure}, but in the context of
the Bayesian network perspective. A naive model for representing
uncertainty in each of the pore-scale parameters assumes pairwise
independent priors. Recall that for
$\B{\Theta} = (\Theta_1, \dots, \Theta_n)$, the variables
$\Theta_1, \dots, \Theta_n$ are pairwise independent,
$\Theta_i \ind \Theta_j$ for all $i,j = 1, \dots, n$ such that
$i \neq j$, if and only if
\begin{equation*}
  \Pb(\Theta_i, \Theta_j) = \Pb(\Theta_i) \Pb(\Theta_j) \,,
\end{equation*}
where $\Pb(\Theta_i)$ denotes the (marginal) PDF of $\Theta_i$ and
$\Pb(\Theta_i, \Theta_j)$ denotes the joint PDF of
$(\Theta_i, \Theta_j)$ (\cite{KollerFriedman:2009gm}). Under the
assumption of independent priors $\Pb(\Theta_i)$, the joint PDF
$\Pb(\B{\Theta})$ factors into the product of the priors,
\begin{equation}
  \label{eq:dist-Theta-ind-params}
  \Pb(\B{\Theta}) = \prod_{i=1}^n \Pb(\Theta_i)\,.
\end{equation}
The full statistical model for $P$ in \cref{eq:statistical-model} is
then given by
\begin{equation}
  \label{eq:P_0}
  P_0 \defeq \delta_{U} (U - u(\B{x},t; \B{X})) \; \delta_{\B{X}}
  (\B{X} - \B{\chi}(\B{\xi},t; \Theta_1, \dots, \Theta_n))
  \prod_{i=1}^n \Pb(\Theta_i)\,.
\end{equation}
In the case of \cref{eq:P_0}, the Bayesian network has the special
form in \cref{fig:general-model-indpriors} where the independent
priors assumption leads to an overall flatness in the graph structure
for $\Pb(\B{\Theta})$. For PDF $\Pb(\Theta_i)$ with finite variance,
the pairwise independence assumption implies the variables $\Theta_i$
are uncorrelated. In particular, \cref{eq:P_0} holds with
$|\B{\Theta}| = n = 4$ for the parameters $\{R, \theta, l, d\}$ that
describe pore-scale features in \cref{fig:pore-structure}. A model or
PDF can then be specified for each $\Theta_i$, for example, the
uniform priors considered in \cite{UmZhangKatsoulakisEtAl:2017aa}.

The Bayesian network in \cref{fig:general-model-indpriors} does not
take into account geometrical structural constraints between the pore
scale parameters that naturally arise when considering a periodic
arrangement of unit cells $\uvol$ as in the hierarchical nanoporous
structure in \cref{fig:pore-structure}. For example, assuming
independent priors $\Pb(\Theta_i)$ each uniformly distributed
according to the hyperparameter ranges in
\cref{tab:comparison-range,tab:extended-range}, it is possible to
sample geometries that are inconsistent with
\cref{fig:pore-structure}. Next we build Bayesian networks based on
causal relationships that encode such natural structural constraints.

\subsection{Correlations arising from pore-scale structural
  constraints}%
\label{sec:correlations-from-constraints}%

\begin{figure}[!h]
  \centering
  \begin{subfigure}[b]{0.3\textwidth}
    \includegraphics[width=1\linewidth]{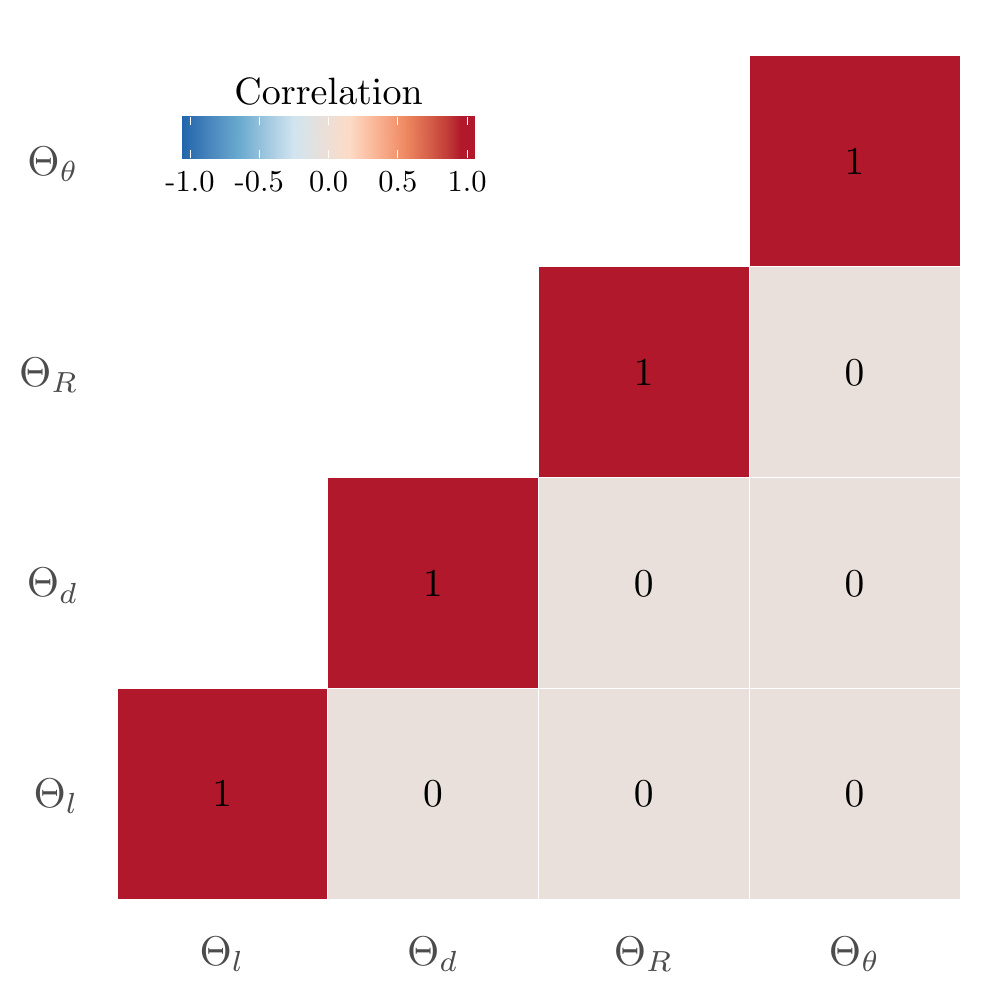}
    \caption{Independent priors \cref{eq:dist-Theta-ind-params} over
      narrow hyperparameter range.}
    \label{fig:corr-exInd-narrow}%
  \end{subfigure}
  \qquad
  \begin{subfigure}[b]{0.3\textwidth}
    \includegraphics[width=1\linewidth]{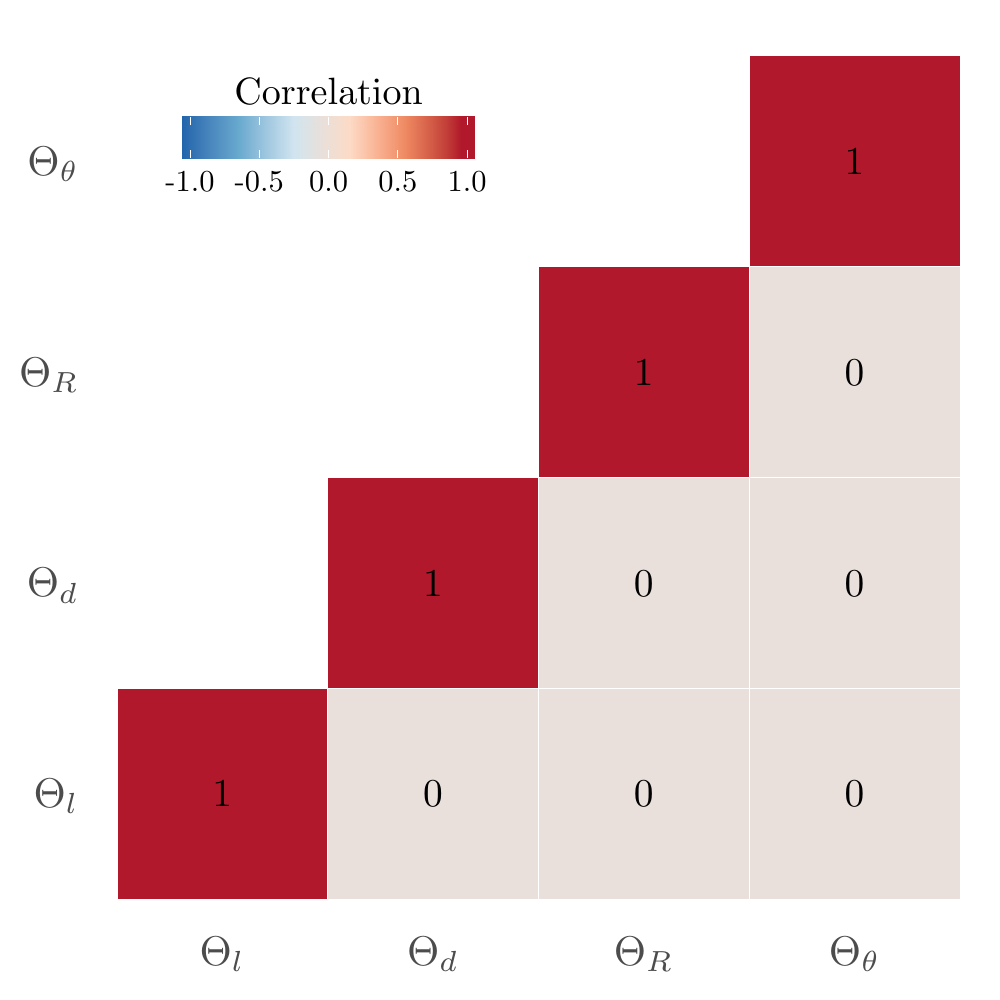}
    \caption{Causal priors \cref{eq:dist-Theta-corr-l-d} over narrow
      hyperparameter range.}
    \label{fig:corr-exB-narrow}%
  \end{subfigure}
  \qquad
  \begin{subfigure}[b]{0.3\textwidth}
    \includegraphics[width=1\linewidth]{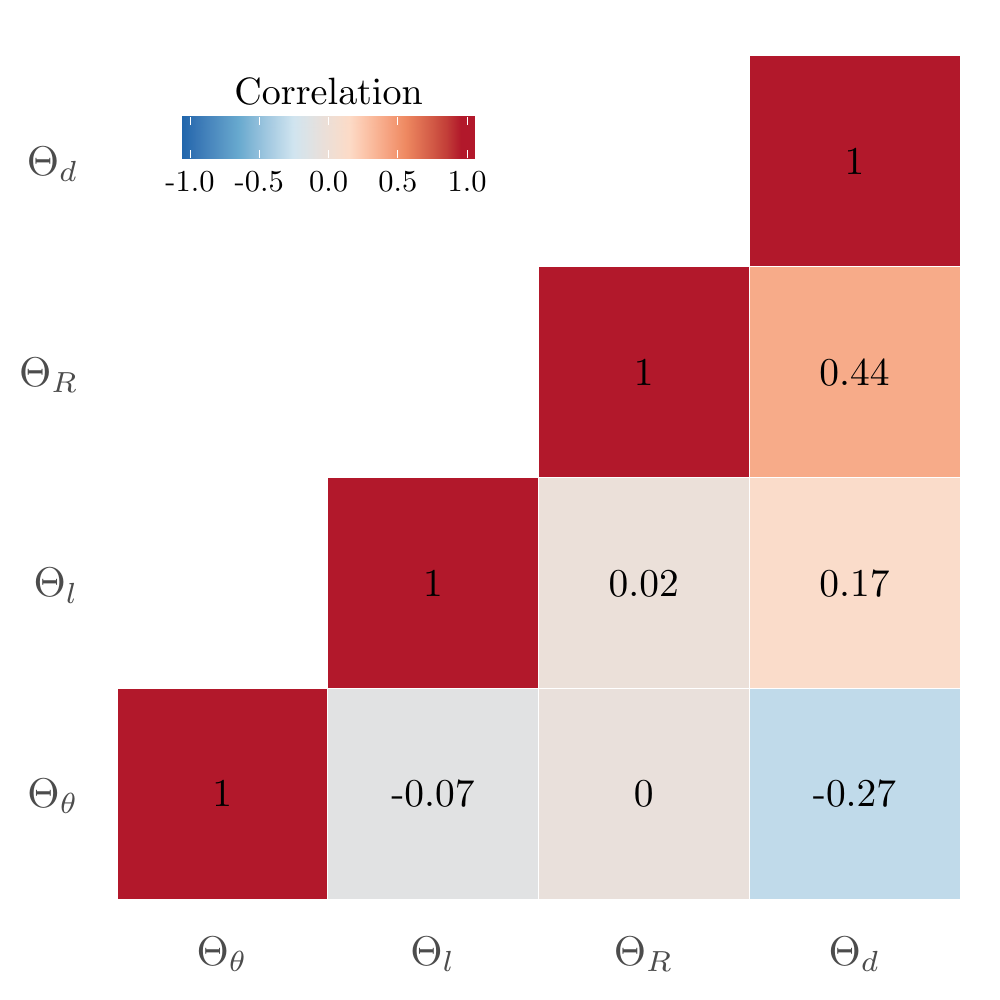}
    \caption{Causal priors \cref{eq:dist-Theta-corr-l-d} over physical
      hyperparameter range.}
    \label{fig:corr-exB-extended}%
  \end{subfigure}
  \qquad
  \caption{The empirical correlation structure of $\Pb(\B{\Theta})$,
    the distribution on pore-scale features, is presented for the
    probabilistic model $P_0$ with independent priors
    \cref{eq:dist-Theta-ind-params} and the model $P_1$ with causal
    priors \cref{eq:dist-Theta-corr-l-d} over both the narrow and
    physical hyperparameter ranges in
    \cref{tab:extended-range,tab:comparison-range}, respectively. In
    general, the physical hyperparameter range is inaccessible to the
    model $P_0$ as the sampling strategy violates structural
    constraints by producing sample geometries inconsistent with
    \cref{fig:pore-structure}. In contrast, the causal relationships
    included in $P_1$ enable sampling over the physical hyperparameter
    range thereby impacting predictions of QoIs (see
    \cref{sec:quant-uncert-darcy}), however the nontrivial correlation
    structure in \cref{fig:corr-exB-extended} poses a challenge for
    global sensitivity analysis (see \cref{sec:global-sa}).}
  \label{fig:correlations}
\end{figure}

Inclusion of causal relationships imposed by geometrical constraints
on the pore-scale parameters adds more complexity to the structure of
$\B{\Theta}$ in the Bayesian network in
\cref{fig:general-model-indpriors}. Here and below we use $\Theta_p$
with labels $p \in \{R, \theta, l, d\}$ in place of indices where no
confusion arises. Moreover, for each $\Theta_p$ we fix hyperparameters
$\{p_{+}, p_{-}\}$ corresponding to upper and lower bounds on
$\Theta_p$ in \cref{tab:comparison-range,tab:extended-range};
performing inference over the Bayesian network to infer
hyperparameters from relevant data is beyond the scope of this work.
In particular, there are a plurality of Bayesian networks that
describe structural constraints consistent with
\cref{fig:pore-structure}. Different causal relationships mirror
distinct engineering or design processes (see
\cref{sec:alternative-prob-models}) and yield different correlation
structures among pore-scale features (cf.\ \cref{fig:correlations} and
\cref{fig:corr-exA-extended}).

\begin{table}[!h]
  \centering
  \caption{Narrow hyperparameter ranges that allow for comparison with
    results in \cite{UmZhangKatsoulakisEtAl:2017aa}.}
  \label{tab:comparison-range}
  \begin{tabular}{lSSSS}
    \toprule
    & \multicolumn{1}{c}{$p = R$ (\si{\nano\meter})} 
    & \multicolumn{1}{c}{$p = \theta$ (\si{\radian})} 
    & \multicolumn{1}{c}{$p = d$ (\si{\nano\meter})} 
    & \multicolumn{1}{c}{$p = l$ (\si{\nano\meter})}\\ 
    \midrule
    (maximum) $p_+$ & 60  & 0.7      & 8   & 18  \\ 
    (minimum) $p_-$ & 10  & 0.07     & 4   & 8   \\
    \bottomrule
  \end{tabular}
\end{table}

\begin{table}[!h]
  \centering
  \caption{Physical hyperparameter ranges.}
  \label{tab:extended-range}
  \begin{tabular}{lSSSS}
    \toprule
    & \multicolumn{1}{c}{$p = R$ (\si{\nano\meter})} 
    & \multicolumn{1}{c}{$p = \theta$ (\si{\radian})} 
    & \multicolumn{1}{c}{$p = d$ (\si{\nano\meter})} 
    & \multicolumn{1}{c}{$p = l$ (\si{\nano\meter})}\\ 
    \midrule
    (maximum) $p_+$ & 60  & 0.4\pi    & 60  & 60  \\ 
    (minimum) $p_-$ & 10  & 0.05\pi   & 5   & 1   \\
    \bottomrule
  \end{tabular}
\end{table}

\begin{figure}[!h]
  \centering %
  \includegraphics[width=.4\textwidth]{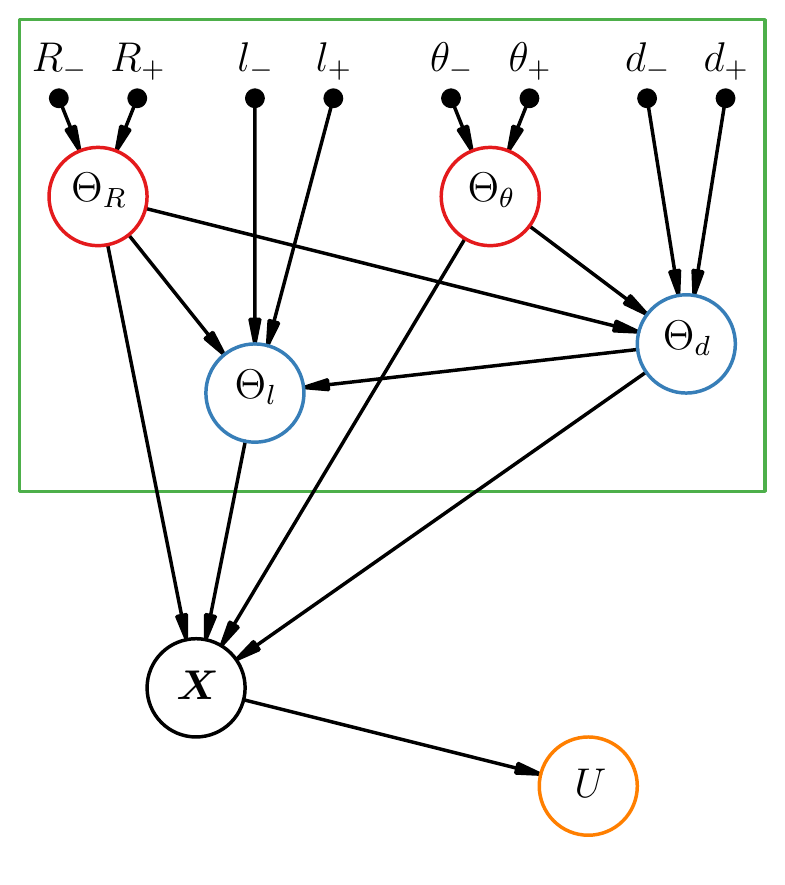}
  \caption{The rich structures of the Bayesian network above,
    representing the probabilistic model $P_1$ in \cref{eq:P_1},
    encodes causal relationships arising from structural constraints
    (cf.\ \cref{fig:upperbound_d}) that are absent in the model $P_0$
    in \cref{eq:P_0} with independent priors in
    \cref{fig:general-model-indpriors}. In this Bayesian network,
    conditional dependencies among the variables $\B{\Theta}$ induce
    various correlation structures that depend on the selected
    hyperparameters (cf.\ the correlation structure for model $P_1$
    over the narrow range of hyperparameters in
    \cref{fig:corr-exB-narrow} vs.\ the physical range in
    \cref{fig:corr-exB-extended}).}
  \label{fig:conditional-d-and-l}
\end{figure}

\begin{figure}[!h]
  \centering %
  \includegraphics[width=.3\textwidth]{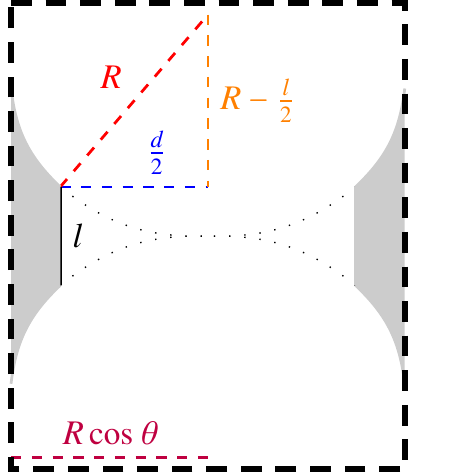}
  \caption{The conditional distribution of $\Theta_d$ and $\Theta_l$
    in \cref{eq:Theta_d-given-R-theta,eq:Theta_l-given-R-theta-d}
    arises from geometric constraints that arise naturally when
    considering the hierarchical nanopore structure in
    \cref{fig:pore-structure}. Above, we illustrate that the nanotube
    radius $d/2$ cannot exceed $R\cos\theta$, half the width of the
    unit cell $\uvol$ resulting in \cref{eq:Theta_d-given-R-theta}.
    Further, to exclude gaps between vertical mesopores that are less
    than zero, $l > 2R - \sqrt{4R^2-d^2}$ from the right triangle in
    the diagram above resulting in \cref{eq:Theta_l-given-R-theta-d}.}
  \label{fig:upperbound_d}
\end{figure}

\Cref{fig:conditional-d-and-l} presents one possible Bayesian network
capturing causal relationships that encode structural constraints
among the hierarchical pore-scale parameters $\{R,\theta,l,d\}$.
Choosing $\{\Theta_R, \Theta_\theta\}$ as independent parameters to be
consistent with our design goal and assuming uniform priors,
\begin{equation}
  \label{eq:Theta_R-Theta_t}
  \Theta_R 
  \sim \unif(R_{-}, R_{+})
  \quad \text{and} \quad 
  \Theta_\theta
  \sim \unif(\theta_{-}, \theta_{+})\,,
\end{equation}
constrains both the distribution of $\Theta_d$ and of $\Theta_l$.
Specifically, in order for the sample $\Theta_d = d$ to be consistent
with the features of the hierarchical nanopore structure in
\cref{fig:pore-structure}, the nanotube diameter cannot exceed the
width of the unit cell $\uvol$, i.e.,\ $d < 2 R \cos \theta$ (see
\cref{fig:upperbound_d}). That is, we are naturally able to specify
the conditional PDF $\Theta_d$ given samples $\Theta_R = R$ and
$\Theta_\theta = \theta$. Assuming an uninformed or uniform model for
this conditional PDF we have
\begin{equation}
  \label{eq:Theta_d-given-R-theta}
  \Theta_d \given \Theta_R, \Theta_\theta \sim 
  \unif(d_-, \min\{2 \Theta_R \cos\Theta_\theta,\, d_+\})\,.
\end{equation}
The length of the nanotube is bounded below by
$l > 2R - \sqrt{4R^2 - d^2}$ when the vertical gap between the
mesopores is zero (see \cref{fig:upperbound_d}) and then the PDF of
$\Theta_l$ given $\Theta_R$ and $\Theta_d$ is
\begin{equation}
  \label{eq:Theta_l-given-R-theta-d}
  \Theta_l \given \Theta_R, \Theta_d \sim 
  \unif(\max\{l_-,\, 2\Theta_R - \sqrt{4\Theta_R^2-\Theta_d^2}\}, l_+)\,,
\end{equation}
where again we assume a uniform model for the conditional PDF of
$\Theta_l$.

The joint distribution for the pore-scale input parameters is then
given by
\begin{equation}
  \label{eq:dist-Theta-corr-l-d}
  \Pb(\B{\Theta}) = \Pb(\Theta_R, \Theta_\theta, \Theta_l, \Theta_d ) 
  = \Pb(\Theta_l \given \Theta_R, \Theta_d) 
  \Pb(\Theta_d \given \Theta_R, \Theta_\theta)
  \Pb(\Theta_R) 
  \Pb(\Theta_\theta) 
\end{equation}
where each of the conditional and prior PDFs is specified in
\eqref{eq:Theta_R-Theta_t}--\eqref{eq:Theta_l-given-R-theta-d} (cf.\
to the assumption of independent priors in
\cref{eq:dist-Theta-ind-params}). Once a range of hyperparameters is
fixed, one can sample from the PDFs
\eqref{eq:Theta_R-Theta_t}--\eqref{eq:Theta_l-given-R-theta-d} and
hence the correlation structure of $\Pb(\B{\Theta})$ can be computed
empirically. In \cref{fig:correlations}, we compare the empirical
correlation structure obtained from \cref{eq:dist-Theta-corr-l-d} over
both a narrow hyperparameter range in \cref{tab:comparison-range} and
a physical hyperparameter range in \cref{tab:extended-range} to the
correlation structure obtained from \cref{eq:dist-Theta-ind-params}.
The model with independent priors is not valid over the physical range
in \cref{tab:extended-range}.

Incorporating \cref{eq:dist-Theta-corr-l-d} into the full statistical
model \cref{eq:statistical-model}, we obtain a new probabilistic model
$P_1$,
\begin{equation}
  \label{eq:P_1}
  P_1 \defeq \delta_{U} (U - u(\B{x},t; \B{X})) \; 
  \delta_{\B{X}} (\B{X} -
  \B{\chi}(\B{\xi},t; \B{\Theta})) \; 
  \Pb(\Theta_l \given \Theta_R, \Theta_\theta, \Theta_d) 
  \Pb(\Theta_d \given \Theta_R, \Theta_\theta)
  \Pb(\Theta_R) 
  \Pb(\Theta_\theta) \,,
\end{equation}
that describes the full statistical model for the multiscale system
with the causal relationships in \cref{fig:conditional-d-and-l}
assuming that uncertainty only enters through the parameters. When the
models for $U$ and $\B{X}$ are known and trivial (e.g., have a known
analytic form), then sampling \cref{eq:P_1} is a straightforward task.
In the next section, we review tools that will enable us to feasibly
sample the statistical models for $U$ and $\B{X}$ in order to carry
out UQ and GSA.

\subsection{Constructing a physics-informed probabilistic model for
  macroscopic QoIs}
\label{sec:unknown-darcy-var-models}
We are interested in functionals
$g = g(\B{\Theta}) = g (\Deff(\B{\Theta}), \geff(\B{\Theta}))$ such as
the projections,
\begin{equation}
  \label{eq:DL-and-DT}
  \DL \defeq \Deff^{11} 
  \qquad \text{and} \qquad 
  \DT \defeq \Deff^{22}\,,
\end{equation}
corresponding to the longitudinal and the transverse components of the
effective diffusion coefficient tensor. In the hierarchical nanoporous
media in \cref{fig:pore-structure}, $\DL$ is related to diffusion
through nanotunnels/mesopores and $\DT$ through nanotubes. In the
numerical experiments presented in the sequel, we report uncertainty
in these macroscopic quantities by giving an estimate of the PDF
$f_g \sim \Pb(g \given \B{\Theta})$ for the PDF of $g$ given knowledge
of pore-scale inputs. In the context of the homogenized pore- to
Darcy-scale model in \cref{fig:general-model}, we have the
representation
\begin{equation}
  \label{eq:g-given-Theta}
  \Pb(g \given \B{\Theta}) = \Pb( g \given \B{X}) 
  \Pb(\B{X} \given \B{\Theta}) \Pb(\B{\Theta})\,.
\end{equation}
For example, the univariate PDF for
$\DL = g(\Deff(\B{\Theta}), \geff(\B{\Theta}))$ based on the Bayesian
network in \cref{fig:conditional-d-and-l} is
\begin{equation}
  \label{eq:dist-DL}
  \Pb(\DL \given \B{\Theta}) 
  =  \Pb(\DL \given \B{X}) \; \delta_{\B{X}} (\B{X} -
  \B{\chi}(\B{\xi},t; \B{\Theta})) \; 
  \Pb(\Theta_l \given \Theta_R, \Theta_\theta, \Theta_d) 
  \Pb(\Theta_d \given \Theta_R, \Theta_\theta)
  \Pb(\Theta_R) 
  \Pb(\Theta_\theta) \,.
\end{equation}
In general, the form of the statistical model or PDF for
$\Pb(g \given \B{X})$ and hence $\Pb(g \given \B{\Theta})$ in
\cref{eq:g-given-Theta} is unknown. However, the PDF in
\cref{eq:g-given-Theta} can be estimated empirically using simulations
of the forward model. Sampling \cref{eq:g-given-Theta} based strictly
on input-output pairs may be computationally expensive due to the high
cost involved in simulating the forward model.

The computation of estimates for QoIs of the form
\cref{eq:g-given-Theta} is made feasible using a two-step method that
relies on first finding a truncated generalized polynomial chaos
expansion (gPCE) for the variable $g$ and second using this surrogate
to build an appropriate kernel density estimator (KDE) for the desired
distribution. A surrogate $\hat{g}$ for $g(\B{\Theta})$ is given by
the gPCE (\cite{LeMaitreKnio:2010sm,Xiu:2010nm,GhanemSpanos:1991pc}),
\begin{equation}
  \label{eq:gPCe-g}
  g(\B{\Theta}) = \sum_{i=0}^\infty G_i \Psi_i (\B{\Theta}) 
  \approx \sum_{i=0}^{N_{PC}} G_i \Psi_i (\B{\Theta}) \eqdef \hat{g}\,,
\end{equation}
where $\Psi_i(\B{\Theta})$ are an orthogonal multivariate polynomial
basis, $G_i$ are the expansion coefficients, and the expansion is
truncated after $N_{PC}$ terms such that
\begin{equation}
  \label{eq:N-PC}
  N_{PC} -1 = \prod_{i=1}^{N_{p}} (1+\kappa_{i}) \,,
\end{equation}
where $\kappa_{i}$ is the polynomial order bound for the
$i^{\mathrm{th}}$ dimension and $N_{p}$ is the number of parameters.

The second step involves producing $N$ samples $\hat{g}^k$ using the
gPCE \cref{eq:gPCe-g} corresponding to realizations $\B{\Theta}^{k}$
for $k = 1, \dots, N$, for a fixed number $N$. These samples are then
used to construct a KDE $\bar{f}_{\hat{g}}$ for the desired density
$f_g$, for example, using a Gaussian-kernel,
\begin{equation}
  \label{eq:general-kde-g}
  \bar{f}_{\hat{g}}(\eta) = \frac{1}{N \sqrt{2\pi h^2}} 
  \sum_{k=1}^{N} \exp\left[-\frac{(\eta-\hat{g}^k)^2}{2h^2}\right]\,,
\end{equation}
where $h$ is the kernel bandwidth; for the multivariate density
$\bar{f}_{\hat{g}_1, \hat{g}_2}$ we similarly employ Gaussian-kernels
\begin{equation}
  \label{eq:general-jkde-g1-g2}
  \bar{f}_{\hat{g}_1, \hat{g}_2} (\eta_1, \eta_2) = 
  \frac{1}{N \sqrt{2 \pi h_1 h_2}} \sum_{k=1}^N 
  \exp \left[ -\frac{(\eta_1 - \hat{g}_1^k)^2}{2h_1^2}
    -\frac{(\eta_2 - \hat{g}_2^k)^2}{2h_2^2} \right]\,,
\end{equation} 
with bandwidths $h_1$ and $h_2$. Software \texttt{DAKOTA}
\cite{DAKOTA:2009man} was used to automatize the process of computing
the coefficients, basis functions, and truncation parameters appearing
in \cref{eq:gPCe-g}. This approach is taken in
\cite{UmZhangKatsoulakisEtAl:2017aa} in the context of independent
priors with the aim of computing Sobol' indices for global sensitivity
analysis. Although generalizations of gPCE that handle correlated
inputs exist (e.g.,\ \cite{NavarroEtAl:2014aa,PaulsonEtAl:2017pc}), we
instead present a recipe for obtaining the desired gPCE in
\cref{eq:gPCe-g} that non-intrusively utilize existing codes and
software packages such as \texttt{DAKOTA} that implement methods
assuming uncorrelated inputs in the next section.

\section{Quantifying uncertainty in Darcy-scale flow variables}
\label{sec:quant-uncert-darcy}

The present section deals with uncertainty quantification for the
Bayesian network PDE model for multiscale porous media outlined in
\cref{sec:pore-to-darcy-model,sec:prob-graph-model}. The causal
relationships encoded by the Bayesian network for the full statistical
model in \cref{fig:general-model} propagate uncertainty from
pore-scale parameters to Darcy-scale variables via the homogenization
map in
\cref{eq:geff,eq:Deff,eq:closure-equation,eq:closure-norm-cond,eq:closure-bdry-cond,eq:eg-closure-cond-chi1,eq:eg-closure-cond-chi2}.
Together, these allow one to study systematically the impact of
microscopic structural uncertainty on macroscopic flow variables. At
present, we incorporate the Bayesian networks constructed in
\cref{sec:prob-graph-model} for informed priors into the random PDE
homogenization framework thereby examining through simulations and
numerical experiments the role of causality in predictions of
Darcy-scale flow variable QoIs. Specifically, we report on numerical
experiments concerning the joint and marginal distributions of
Darcy-scale flow variables where uncertainty stems from pore-scale
features with correlations arising from structural constraints encoded
by the Bayesian network in \cref{fig:conditional-d-and-l}. As these
numerical experiments utilize \texttt{DAKOTA} (\cite{DAKOTA:2009man})
to compute gPCE surrogates for Darcy-scale variables, we first
illustrate in \cref{sec:rosenblatt} a technique for decorrelating
inputs to allow the non-intrusive use of existing codes and software
packages. This work-flow, which relies on Rosenblatt transformations,
gPCEs, and \texttt{DAKOTA}, enables uncertainty quantification and
global sensitivity analysis by making it feasible to sample from the
desired QoI with respect to a given statistical model $P$ in
\cref{eq:statistical-model}.

\subsection{Non-intrusive input decorrelation using Rosenblatt
  transforms}
\label{sec:rosenblatt}

Many variance-based methods for uncertainty quantification and
sensitivity analysis, such as Sobol' indices, and hence the popular
software packages that implement these methods, require models that
assume statistically independent inputs. Bayesian networks, recall
\cref{fig:general-model-indpriors,fig:conditional-d-and-l}, encode
correlations through the specification of causal relationships (see
\cref{fig:correlations}). Presently, we highlight how the Rosenblatt
transform can be used to decorrelate inputs by mapping a vector of
random variables with a specified joint distribution onto a vector of
independent uniform random variables when the conditional
distributions are known. This procedure, in \cref{alg:correlate}
below, enables the non-intrusive use of \texttt{DAKOTA} and existing
codes for solving the forward model for our application of interest
when the conditional dependencies are represented using Bayesian
networks.

The Rosenblatt transform (\cite{Rosenblatt:1952rt}) turns the problem
of sampling a general joint distribution into the problem of sampling
a vector of independent $\unif(0,1)$ random variables. Let
$\B{X} = (X_1, \dots, X_k)$ be a random vector with a continuous joint
cumulative distribution function $F(x_1, \dots, x_k)$. Define a
transform,
$\mathcal{T}(x) = \mathcal{T}(x_1, \dots x_k) = (z_1, \dots, z_k) =
z$, given by
\begin{equation}
  \label{eq:Rosenblatt}
  \begin{split}
    z_1 &= \Pb(X_1 \leq x_1) = F_1(x_1)\,,\\
    z_2 &= \Pb(X_2 \leq x_2 \given X_1 = x_1)
    = F_{2\given 1} (x_2 \given x_1)\,,\\
    \vdots & \\
    z_k &= \Pb(X_k \leq x_k \given X_{k-1} = x_{k-1}, \dots, X_1 =
    x_1) = F_{k \given k-1, \dots, 1}(x_k \given x_{k-1}, \dots,
    x_1)\,,
  \end{split}
\end{equation}
where $F_{i\given j}$ is the conditional cumulative distribution
function of $X_i$ given $X_j$, i.e.\
$F_{i\given j} (x_i | x_j) = \Pb(X_i < x_i \given X_j = x_j)$. The
Rosenblatt transform, $\B{Z} \defeq \mathcal{T}(\B{X})$, yields
$\B{Z} = (Z_1, \dots, Z_k)$ uniformly distributed on the
$k$-dimensional hypercube, that is, $Z_1, \dots, Z_k$ are independent
and identically distributed (\iid) $\unif(0,1)$ random variables. Note
that this transform depends on the ordering of the elements in the
vector $\B{X}$ due to the serial nature of the conditioning; we denote
the Rosenblatt transform and the inverse, when it exists, associated
with the ordering of a particular vector $\B{X}$ with a subscript,
e.g.\ $\mathcal{T}_{\B{X}}$.

For our application of interest, a target vector $\B{\Theta}$ with a
causal structure encoded by a Bayesian network can be obtained by
applying the inverse Rosenblatt transform to a vector of independent
uniform variables. For example, given the random vector of parameters
$\B{\Theta} = (\Theta_R,\Theta_\theta,\Theta_l,\Theta_d)$ with joint
distribution in \cref{eq:dist-Theta-corr-l-d}, the transform
\cref{eq:Rosenblatt} simplifies to
\begin{align*}
  z_1 &= F_R(x_1) \,,\\
  z_2 &= F_{\theta \given R} (x_2 \given x_1) = F_{\theta}(x_2) \,,
        \qquad (\text{since } \Theta_R \ind \Theta_\theta)\,, \\
  z_3 &= F_{d \given R,\theta}(x_3 \given x_2, x_1) \,,\\
  z_4 &= F_{l \given d, R}(x_4 \given x_3, x_1) \,,
        \qquad \left(\text{since } 
        \Theta_l \ind \Theta_\theta \given \Theta_d \right)\,,
\end{align*}
due to the independence and conditional independence of the variables
(\cite{KollerFriedman:2009gm}), as indicated in
\cref{fig:conditional-d-and-l} by the absence of a causal
relationships between $\Theta_R$ and $\Theta_\theta$ and between
$\Theta_l$ and $\Theta_d$, respectively. Thus, the Rosenblatt
transform is given by
\begin{equation*}
  \mathcal{T}_{\B{\Theta}}(\B{\Theta}) = \mathcal{T}_{\B{\Theta}}(\Theta_R, \Theta_\theta, 
  \Theta_l, \Theta_d) 
  = \left(F_R(\Theta_R), \; F_{\theta}(\Theta_\theta), \;
    F_{d \given R,\theta} (\Theta_d \given \Theta_R, \Theta_\theta), \; F_{l \given R,d}(\Theta_l
    \given \Theta_R, \Theta_d)\right) \eqdef \B{Z}
\end{equation*}
where, using the statistical models indicated in
\cref{eq:Theta_R-Theta_t,eq:Theta_d-given-R-theta,eq:Theta_l-given-R-theta-d},
the components of $\B{Z}$ are
\begin{align*}
  Z_1 &= F_R(\Theta_R) 
        = \frac{\Theta_R - R_{-}}{R_{+} - R_{-}}\,,\\
  Z_2 &= F_\theta (\Theta_\theta) 
        = \frac{\Theta_\theta - \theta_{-}}{\theta_{+} - \theta_{-}}\,,\\
  Z_3 &= F_{d}(\Theta_d \given \Theta_R, \Theta_\theta) 
        = \frac{\Theta_d - d_{-} }{\min\{\Theta_R \cos \Theta_\theta , d_{+}\} - d_{-}}\,,\\
  Z_4 &= F_{l \given R, d}(\Theta_l \given \Theta_R, \Theta_d) 
        = \frac{\Theta_l 
        - \max\{l_-, 2\Theta_R - \sqrt{4\Theta_R^2-\Theta_d^2}\}}
        {l_{+} - \max\{l_-, 2\Theta_R - \sqrt{4\Theta_R^2-\Theta_d^2}\}}\,.
\end{align*}
The corresponding inverse Rosenblatt transform
$\B{\Theta} = \mathcal{T}_{\B{\Theta}}^{-1}(\B{Z})$ is,
component-wise,
\begin{align*}
  \Theta_R &= Z_1(R_{+} - R_{-}) + R_{-} \\
  \Theta_\theta &= Z_2(\theta_{+} - \theta_{-}) + \theta_{-} \\
  \Theta_d &= Z_3( \min\{2\Theta_R \cos\Theta_\theta,\, d_+\} -
             d_{-}) + d_{-}\\
  \Theta_l &= Z_4(l_+ - \max\{l_-,\, 2\Theta_R -
             \sqrt{4\Theta_R^2-\Theta_d^2}\}) + \max\{l_-,\, 2\Theta_R - \sqrt{4\Theta_R^2-\Theta_d^2}\}\,,
\end{align*}
where the target vector $\B{\Theta}$ has the distribution
$\Pb(\B{\Theta})$ given in \cref{eq:dist-Theta-corr-l-d}. Thus, the
inverse transform maps a random vector $\B{Z} \sim \unif(0,1)^k$ into
a target distribution $\Pb(\B{\Theta})$ using knowledge of the
conditional dependencies associated with the Bayesian network for
$\Pb(\B{\Theta})$, in particular, where the conditional distribution
$\Theta_j$ given $\Theta_{j-1}, \dots, \Theta_1$ is known for each
$j = 1, \dots, k$.

Together, Bayesian networks and Rosenblatt transforms provide a
strategy for non-intrusively incorporating constraints and
correlations into existing computational frameworks.
\Cref{alg:correlate} describes the use of \texttt{DAKOTA} for
computing surrogates for Darcy-scale QoIs based on correlated inputs
given an inverse Rosenblatt transform and an existing solver for the
forward problem. A surrogate $\hat{g}$ for a QoI $g(Y(\B{\Theta}))$,
e.g.\ the gPCE coefficients $G$ in \cref{eq:gPCe-g} and truncation
parameter $N_{PC}$ in \cref{eq:N-PC}, can be computed with
\texttt{DAKOTA} using several forward simulations of the response
$Y = \mathcal{M}(\B{\Theta})$ where $\mathcal{M}$ denotes the portion
of the forward model solver that maps a random input $\B{\Theta}$ to
$Y$. For example, if the QoI is $g(Y(\B{\Theta})) = \DL$ then
$\mathcal{M}$ would correspond to the projection $\Deff^{11}$ in
\cref{eq:DL-and-DT} of the solution to the coupled system
\cref{eq:Deff} with
\cref{eq:closure-equation,eq:closure-norm-cond,eq:closure-bdry-cond}.
The compositional model
$Y = (\mathcal{M} \circ \mathcal{T}^{-1})(\B{Z})$ that first employs
the inverse Rosenblatt transform provides a non-intrusive means of
computing with \texttt{DAKOTA} since the statistics of the output of
$\mathcal{M}\circ \mathcal{T}^{-1}$ in response to $\B{Z}$ are
identical to the statistics of the output of $\mathcal{M}$ in response
to $\B{\Theta}$ (\cite{TorreEtAl:2017aa}). An important observation is
that \cref{alg:correlate} returns a surrogate $\hat{g}$ for
$g(Y(\B{\Theta}))$ with respect to the input variables $\B{Z}$, not
for $\B{\Theta}$ directly, and therefore care must be taken if
reporting Sobol' indices with respect to the input parameters
(\cite{MaraEtAl:2015np}). In the next section, we observe that the
density functions for effective Darcy-scale QoIs exhibit non-Gaussian
behavior. Together, these two challenges---correlated inputs and
non-Gaussian QoIs---motivate the investigation of moment-independent
global sensitivity indices in \cref{sec:global-sa}.

\begin{algorithm}[!h]
  \DontPrintSemicolon
  \SetKwInOut{Input}{input}\SetKwInOut{Output}{output}
  \SetKwComment{Comment}{$\triangleright$\ }{}%

  \Input{$\mathcal{T}_{\B{\Theta}}^{-1}$ \Comment*{inverse Rosenblatt
      transform}}%
  \Input{$\mathcal{M}$ \Comment*{forward model solver}}%
  \Output{$\hat{g}_M$ \Comment*{surrogate for $g(Y(\B{\Theta}))$ e.g.\
      gPCE in \cref{eq:gPCe-g}}}%
  \BlankLine %
  \Begin{%
    \emph{\texttt{DAKOTA} as wrapper to produce surrogate using $M$
      input-output simulations}\;%
    \For{$i \leftarrow 1$ \KwTo $M$}{ sample \iid\,
      $Z_k \sim \unif(0,1)$, $k=1,\dots,n$\;
      $\B{Z}_i \leftarrow (Z_1, \dots, Z_n)$\;%
      $\B{\Theta} \leftarrow \mathcal{T}_{\B{\Theta}}^{-1} (\B{Z}_i)$\;%
      $Y_i \leftarrow \mathcal{M}(\B{\Theta})$
      \Comment*{\texttt{DAKOTA} maps independent $\B{Z}_i \mapsto Y_i$}%
      
    }
    $\hat{g}_M \leftarrow \texttt{DAKOTA}\left(Y_1(\B{Z}_1), \dots,
      Y_M(\B{Z}_M)\right)$\;%
    \Return{$\hat{g}_M$} \Comment*{based on $M$ input-output
      simulations}%
  }%
  \caption{Decorrelate inputs via Rosenblatt transforms for
    non-intrusive scientific computing}%
  \label{alg:correlate}
\end{algorithm}

\subsection{Numerical experiments: incorporating causality in
  predictions of Darcy-scale flow}
\label{sec:Darcy-scale-qoi}

The numerical experiments that follow employ the following common
setup that makes sampling the distribution of QoIs computationally
feasible for uncertainty quantification and global sensitivity
analysis. \texttt{DAKOTA} is used as a wrapper to map random inputs on
pore-scale parameters to Darcy-scale responses yielding gPCE
surrogates for Darcy-scale variables, i.e.\ effective longitudinal
diffusion $\DL$, effective transverse diffusion $\DT$, and effective
sorption rate constant $\geff$.

With regard to sampling input-output pairs for the generation of gPCE
surrogates, we follow the decorrelation procedure described in
\cref{alg:correlate} for each sample $\B{\Theta}$ that relates to a
possible configuration of pore-scale features consistent with the
hierarchical nanoporous material in \cref{fig:pore-structure}. In
particular, \cref{alg:correlate} allows seamless, non-intrusive
integration with existing codes for the numerical solution of the
multiscale forward model presented in \cref{sec:pore-to-darcy-model}.
For the numerical solutions of the forward model, we first solve the
closure equations
\cref{eq:closure-equation,eq:closure-norm-cond,eq:closure-bdry-cond,eq:eg-closure-cond-chi1,eq:eg-closure-cond-chi2}
using a finite element code written in \texttt{COMSOL} and then
compute the rate constant $\geff$ in \cref{eq:geff} and effective
diffusions $\DL$ and $\DT$ in \cref{eq:DL-and-DT} by numerically
evaluating the quadrature in \cref{eq:Deff}. For the required gPCE
surrogates, given in \cref{eq:gPCe-g}, we select for $\Psi_i$ the
Askey scheme of hypergeometric orthogonal polynomials
(\cite{Xiu:2010nm}) with $N_{PC} = \num{626}$ (i.e., for $N_p = 4$
parameters we consider polynomials of degree $\kappa_{i}=4$, for each
$i = 1, \cdots, N_{p}$, in \cref{eq:N-PC}). These gPCE surrogates are
then used to construct KDEs for the desired Darcy-scale flow
variables. The KDE approximations for densities below employ the
Gaussian-kernels for univariate and multivariate densities, described
in \cref{eq:general-kde-g,eq:general-jkde-g1-g2}, respectively, with
$N = \num{e8}$ where the kernel bandwidths are estimated using a
modified Sheather-Jones method (\cite{BaringhausFranz:2004tt}).

\begin{figure}[!h]
  \centering
  \begin{subfigure}[t]{0.3\textwidth}
    \includegraphics[width=1\linewidth]{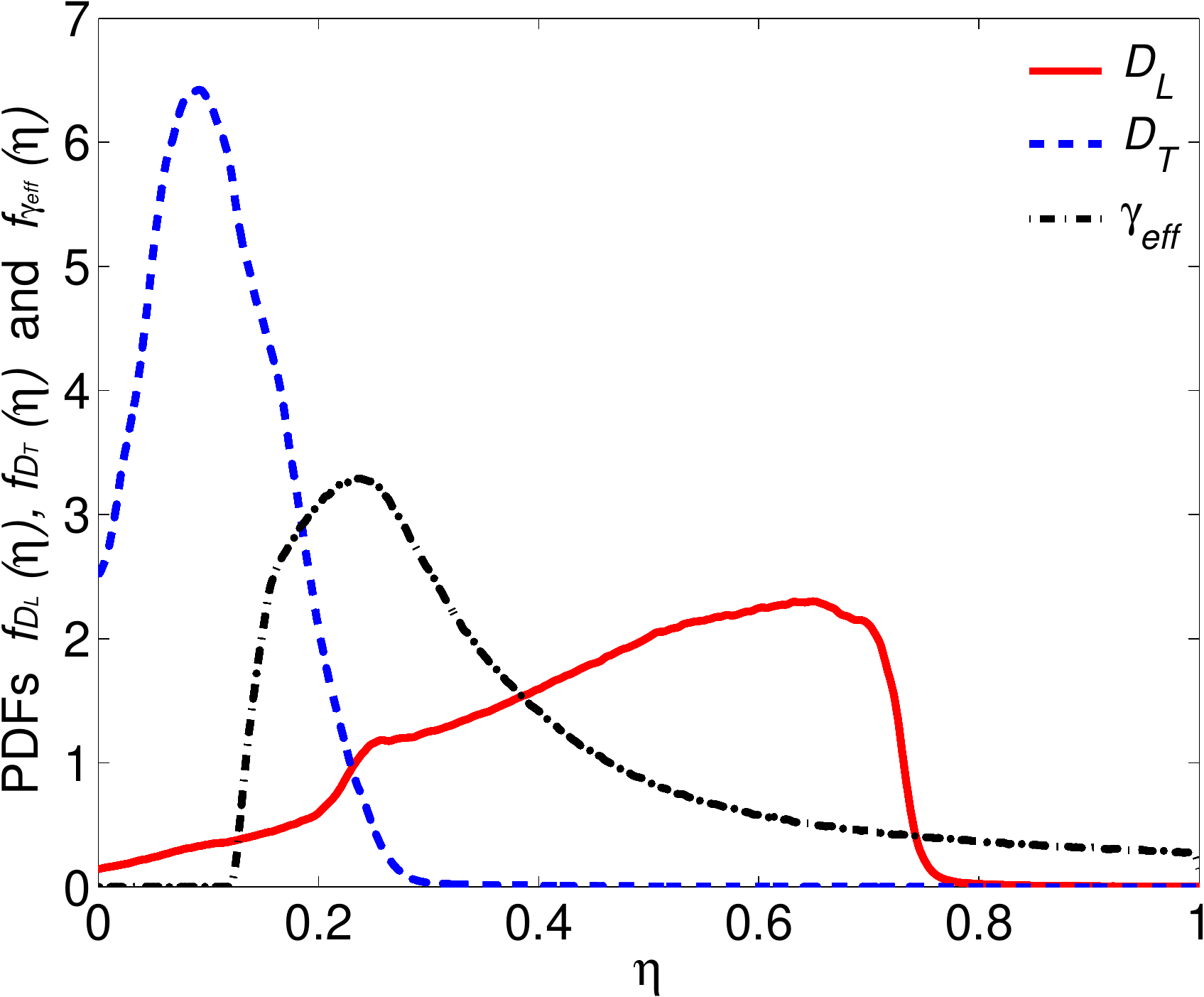}
    \caption{Independent priors (model $P_0$) over narrow
      hyperparameter range.}
    \label{fig:comparison-kdes-indep}%
  \end{subfigure}
  \qquad
  \begin{subfigure}[t]{0.3\textwidth}
    \includegraphics[width=1\linewidth]{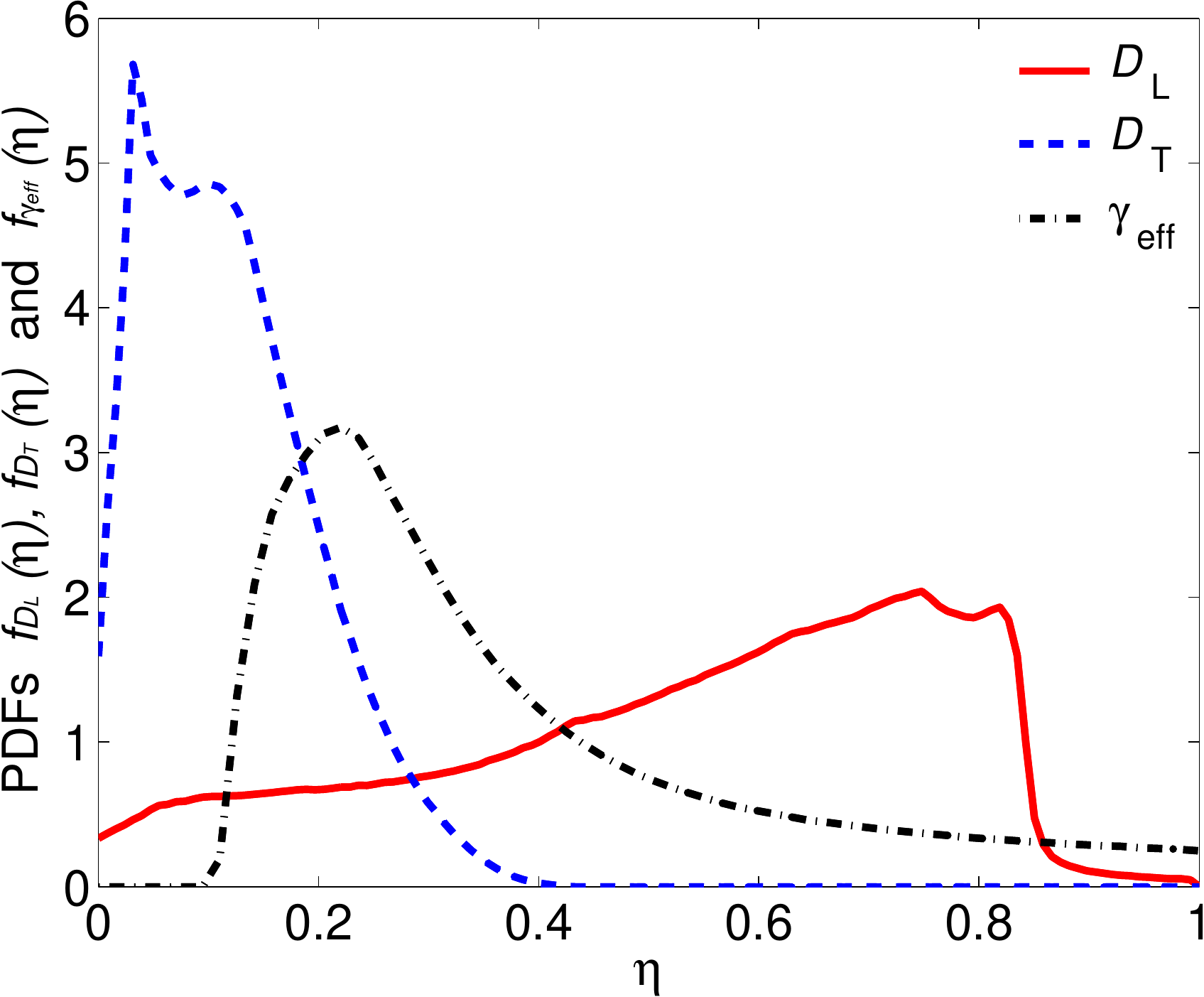}
    \caption{Causal priors (model $P_1$) over narrow hyperparameter
      range.}
    \label{fig:comparison-kdes-corr-exB}%
  \end{subfigure}
  \qquad
  \begin{subfigure}[t]{0.3\textwidth}
    \includegraphics[width=1\linewidth]{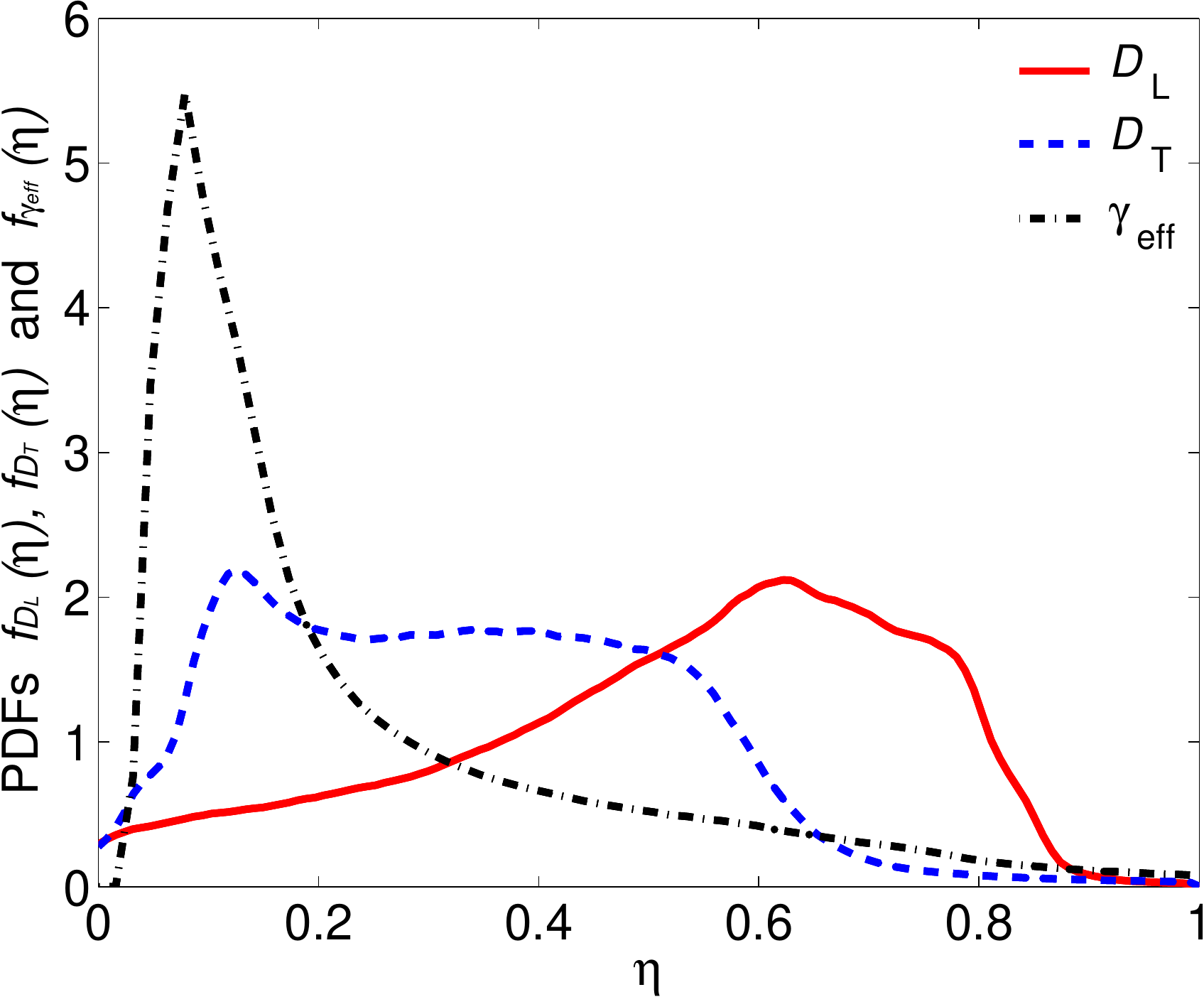}
    \caption{Causal priors (model $P_1$) over physical hyperparameter
      range.}
    \label{fig:comparison-marginal-kdes-exBext}
  \end{subfigure}
  \caption{A comparison of the marginal densities for Darcy-scale QoIs
    above highlights the importance of incorporating causal
    relationships into the modeling process; the distribution of
    Darcy-scale QoIs in \cref{fig:comparison-marginal-kdes-exBext} for
    model $P_1$ in \cref{eq:P_1} with causal priors over the physical
    hyperparamter range in \cref{tab:extended-range} are markedly
    different from the QoIs in \cref{fig:comparison-kdes-indep} for
    model $P_2$ in \cref{eq:P_0} with independent priors over the
    narrow hyperparamter range in \cref{tab:comparison-range}. The
    QoIs in \cref{fig:comparison-kdes-corr-exB} for model $P_1$ are
    expected to be qualitatively similar to the QoIs in
    \cref{fig:comparison-kdes-indep} due to the similarities in the
    correlation structure for the priors over the narrow range of
    hyperparameters (cf.\ \cref{fig:correlations}).}
  \label{fig:comparison-kdes}
\end{figure}

\begin{figure}[!h]
  \centering
  \begin{subfigure}[b]{1\textwidth}
    \includegraphics[width=1\linewidth]{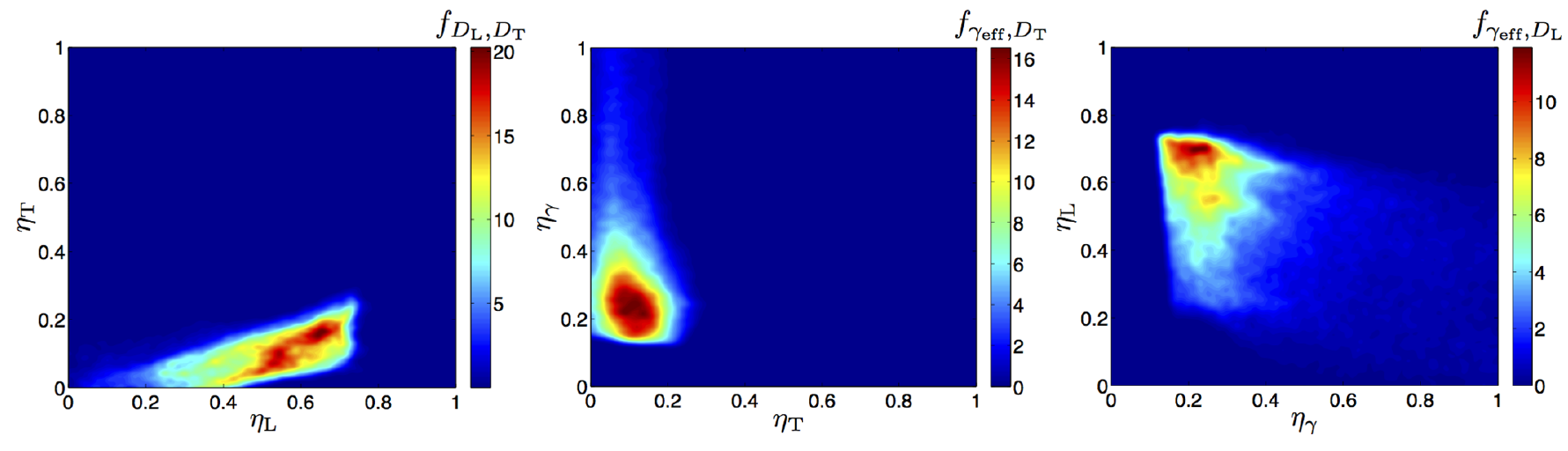}
    \caption{Independent priors (model $P_0$) over narrow
      hyperparameter range, from
      \cite{UmZhangKatsoulakisEtAl:2017aa}.}
    \label{fig:comparison-jkdes-indep}%
  \end{subfigure}\vspace*{2ex}
  
  \begin{subfigure}[b]{1\textwidth}
    \includegraphics[width=.33\linewidth]{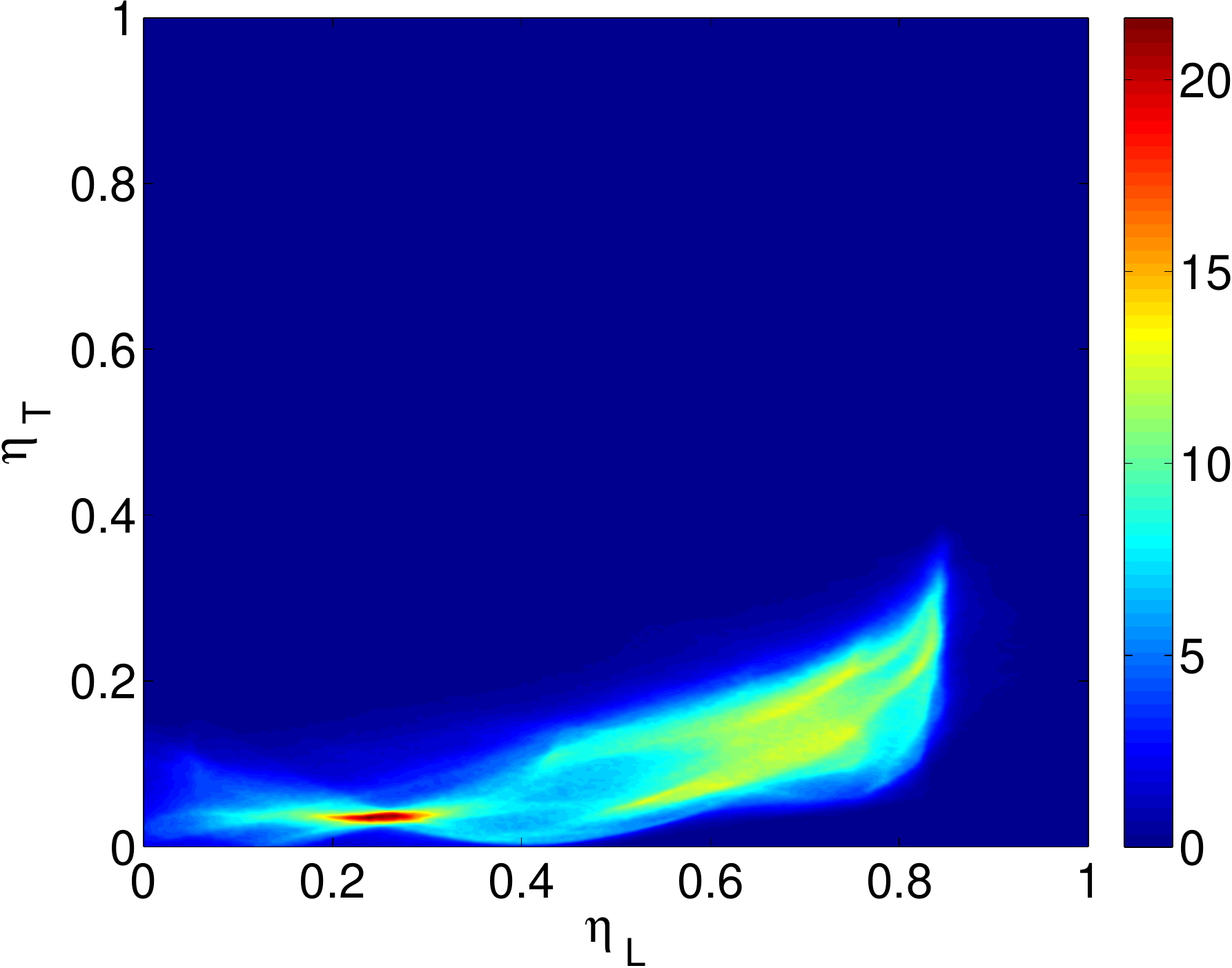}
    \includegraphics[width=.33\linewidth]{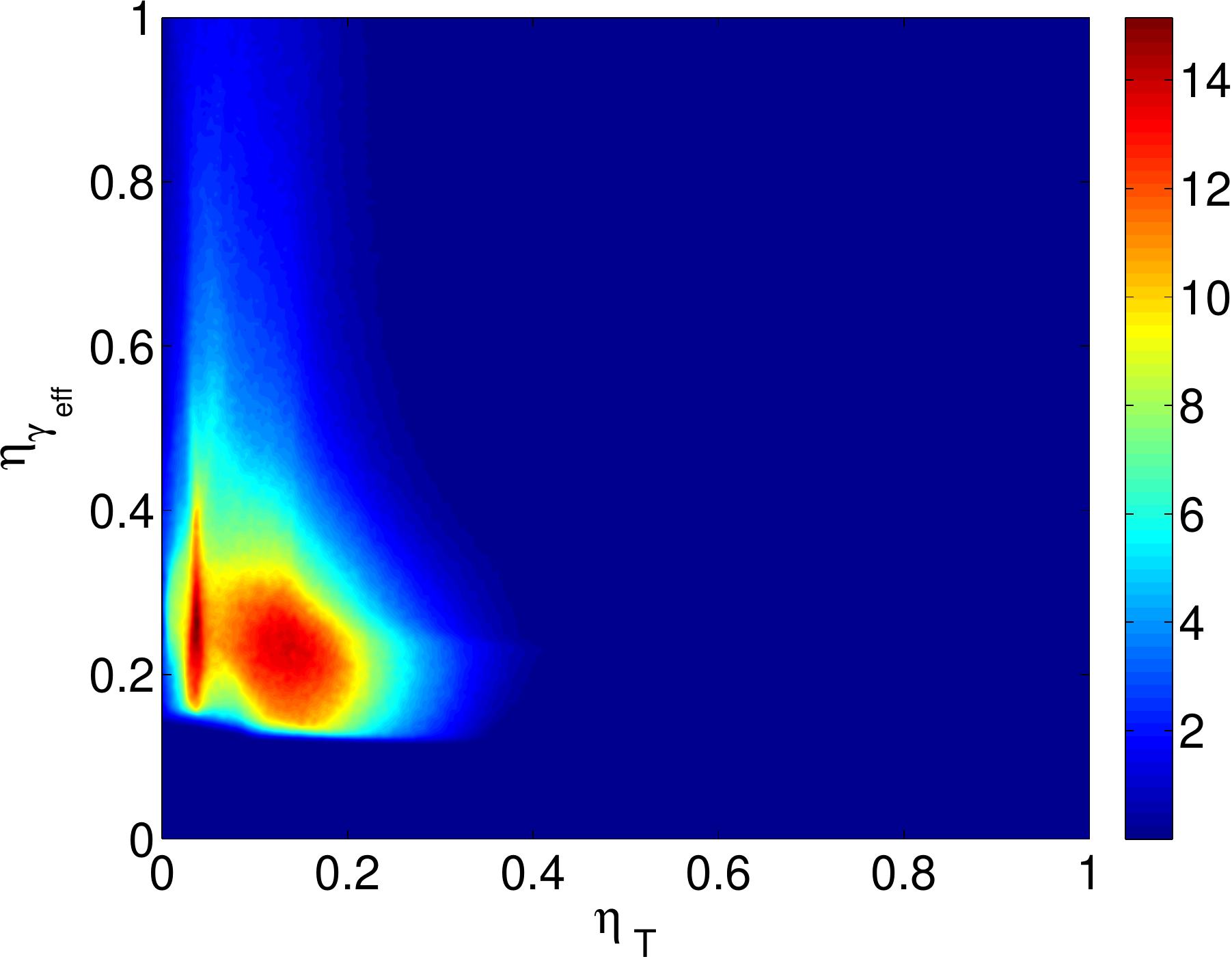}
    \includegraphics[width=.33\linewidth]{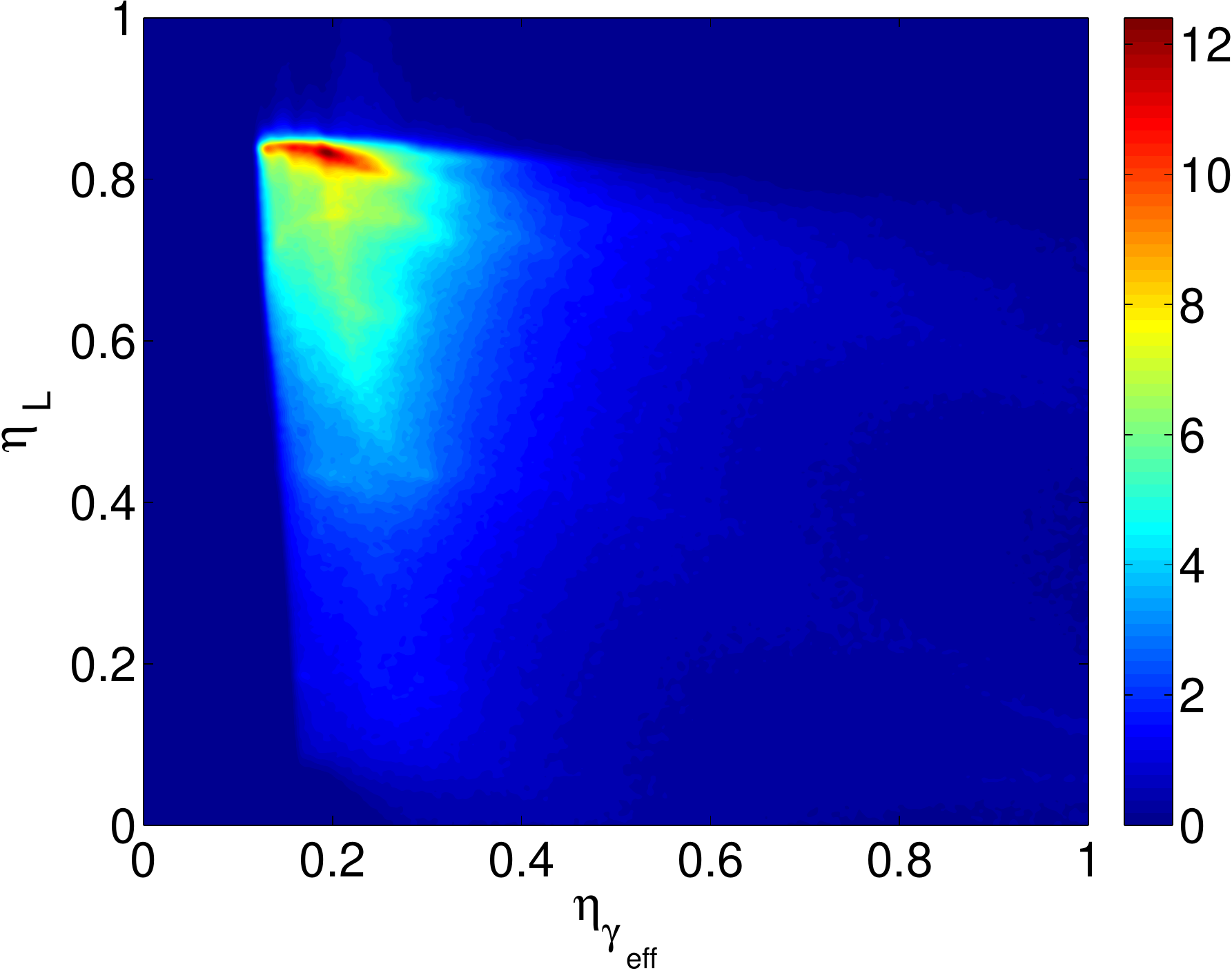}
    \caption{Causal priors (model $P_1$) over narrow hyperparameter
      range.}
    \label{fig:comparison-jkdes-corr-exB}%
  \end{subfigure}\vspace*{2ex}

  \begin{subfigure}[b]{1\textwidth}
    \includegraphics[width=.33\linewidth]{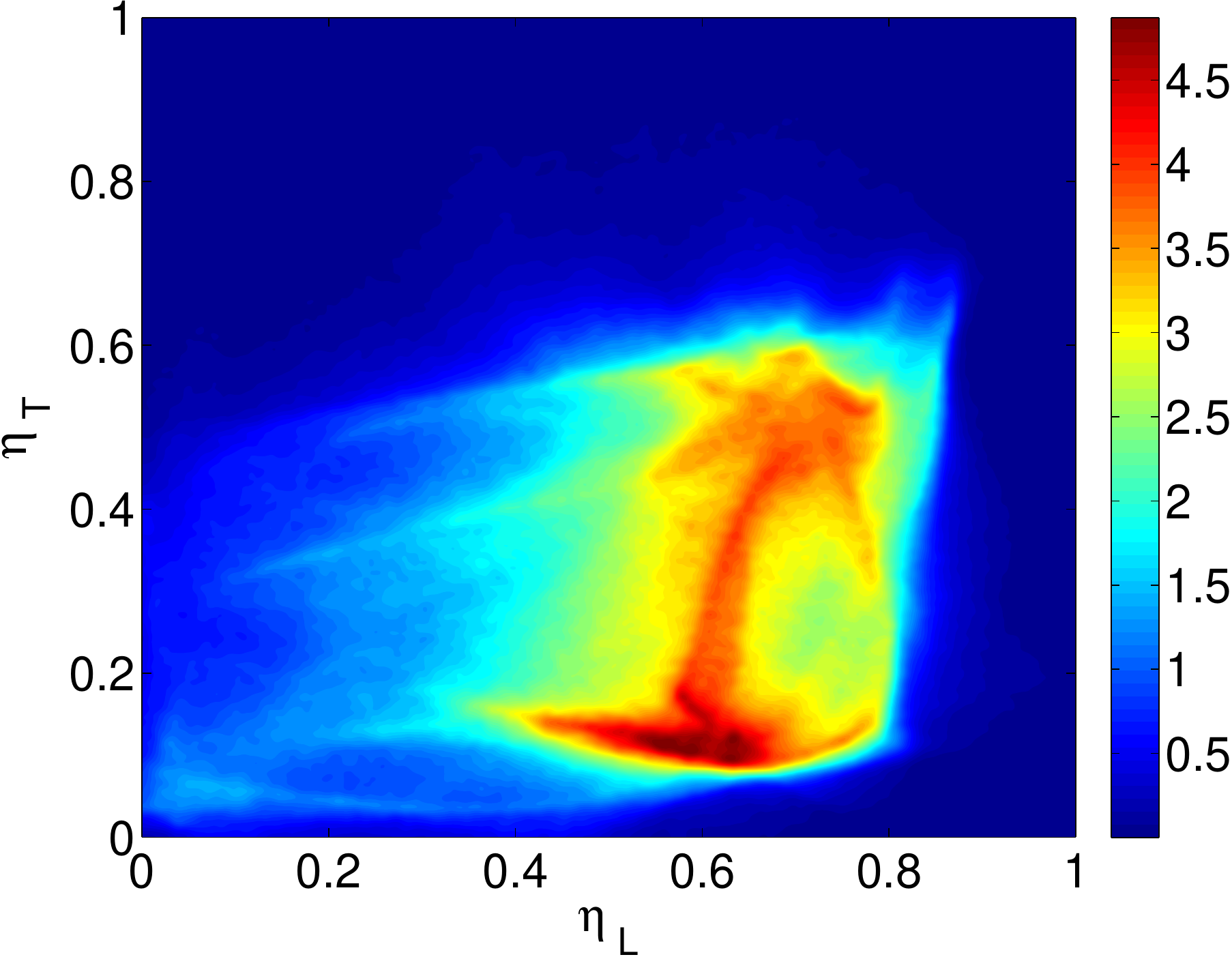}
    \includegraphics[width=.33\linewidth]{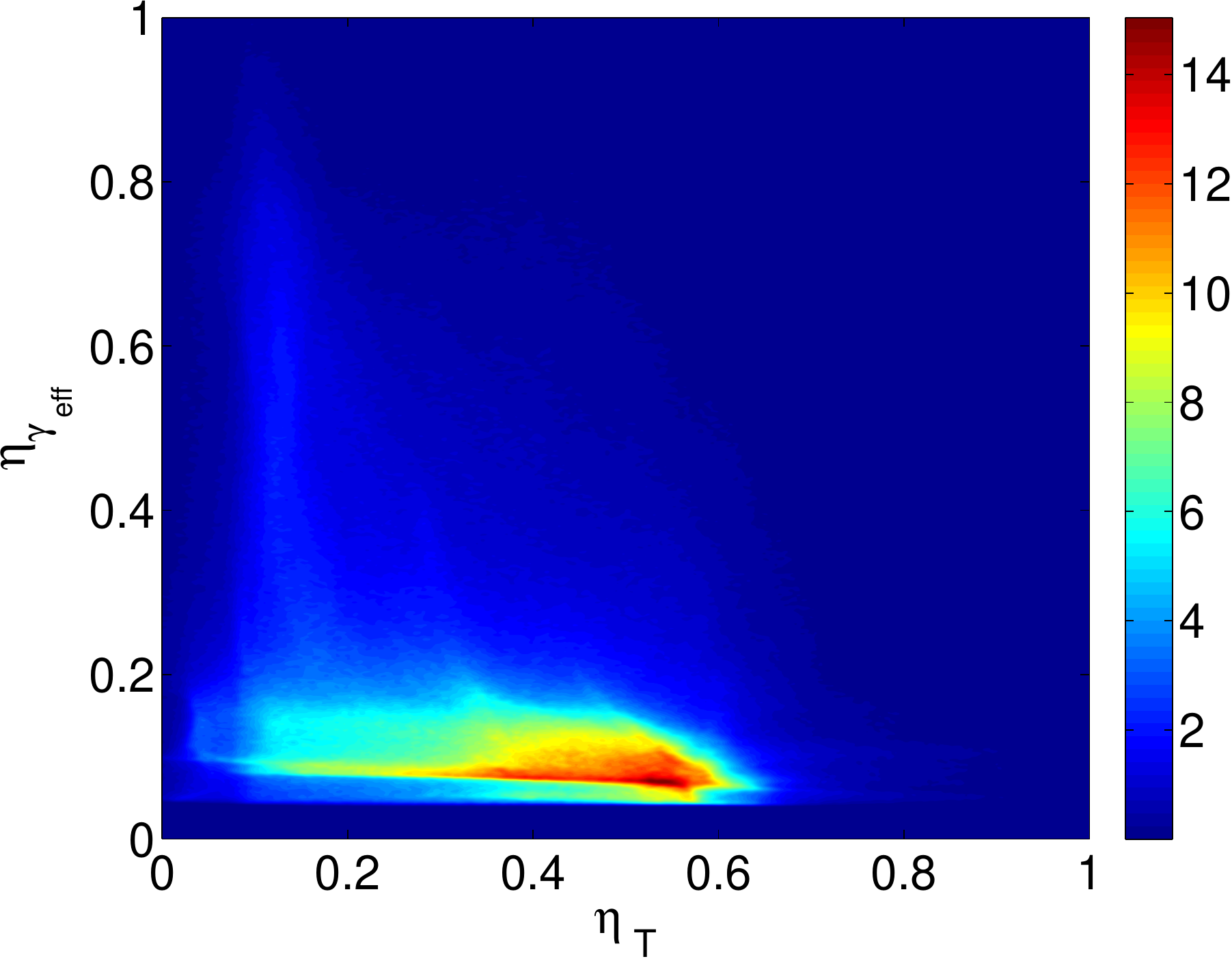}
    \includegraphics[width=.33\linewidth]{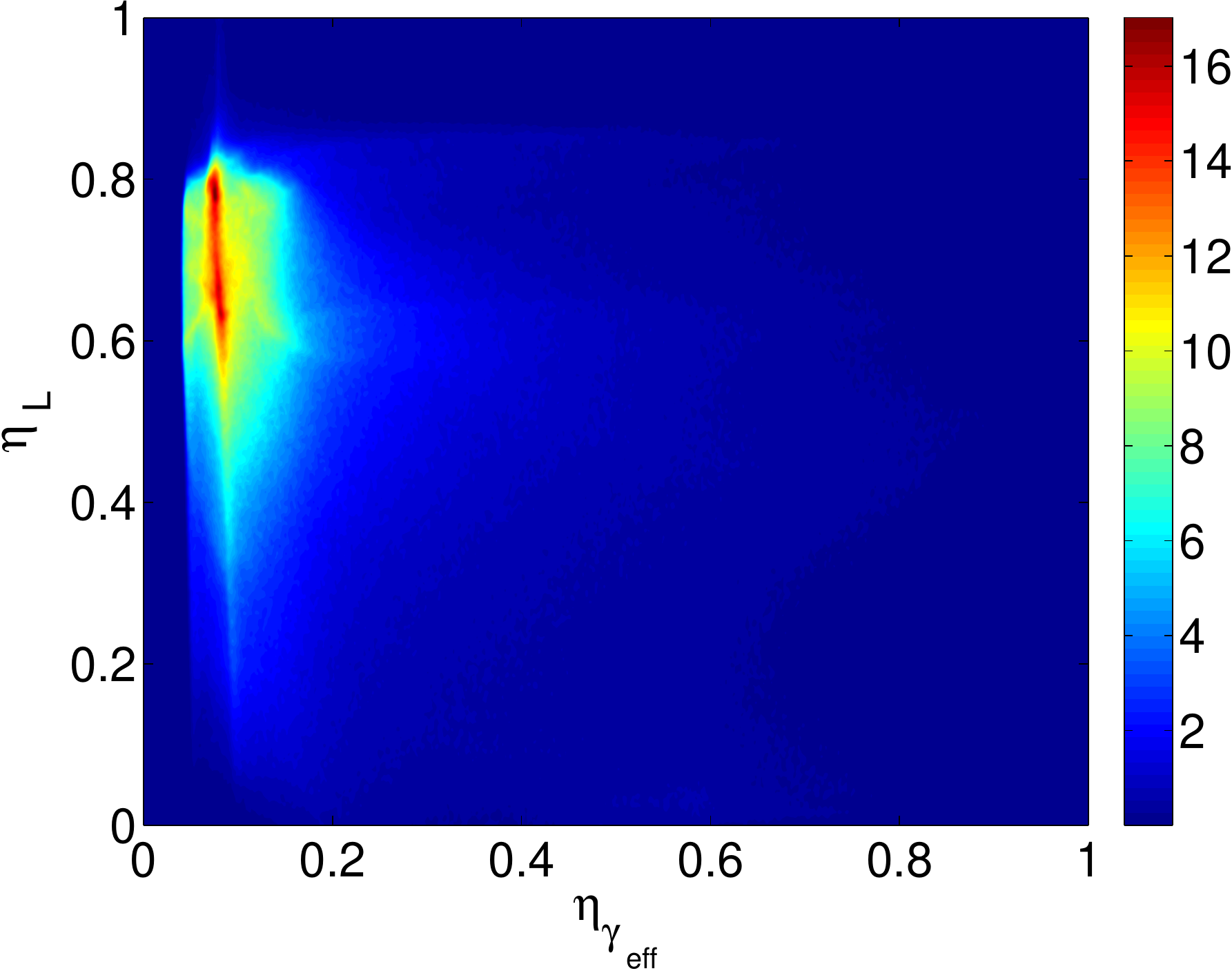}
    \caption{Causal priors (model $P_1$) over physical hyperparameter
      range.}
    \label{fig:comparison-jkdes-corr-exBext}%
  \end{subfigure}

  \caption{A comparison of the joint densities for Darcy-scale flow
    variables above further underscores the importance of including
    causal relationships in the modeling process due to the impact for
    decision support. The causal relationships included in model $P_1$
    in \cref{eq:P_1} guarantee that the pore-scale geometries sampled
    under $P_1$ are consistent with the hierarchical nanoporous
    material in \cref{fig:pore-structure} over the physical
    hyperparamter range in \cref{tab:extended-range}. In
    \cref{fig:comparison-jkdes-corr-exBext}, we observe that the QoIs
    related to model $P_1$ with physical hyperparamters realize a
    richer range of transverse diffusions than the QoIs depicted in
    \cref{fig:comparison-jkdes-corr-exB,fig:comparison-jkdes-indep},
    which correspond to the models $P_0$ and $P_1$ over the narrow
    hyperparameter range in \cref{tab:comparison-range}, thereby
    differentially impacting decision tasks.}
  \label{fig:comparison-joint-kdes}
\end{figure}

As a first experiment, we compare simulations of Darcy-scale flow
variables based on the causal relationships encoded by model $P_1$ in
\cref{eq:P_1} to simulations based on the independent priors model
$P_0$ in \cref{eq:P_0} (cf.\ Bayesian networks in
\cref{fig:conditional-d-and-l,fig:general-model-indpriors}). Over the
narrow range of hyperparameter values given in
\cref{tab:comparison-range}, we anticipate qualitative similarities in
the resulting Darcy-scale outputs as both $P_0$ and $P_1$ exhibit
statistically uncorrelated parameter distributions over this
hyperparameter range as demonstrated by the empirical correlations in
\cref{fig:corr-exInd-narrow,fig:corr-exB-narrow}. In
\cref{fig:comparison-kdes}, we observe that the marginal distributions
for Darcy-scale flow variables based on model $P_1$ in
\cref{fig:comparison-kdes-corr-exB} are qualitatively consistent with
the marginals based on model $P_0$ in
\cref{fig:comparison-kdes-indep}. Likewise, in
\cref{fig:comparison-joint-kdes} the joint distributions $(\DL, \DT)$,
$(\DT, \geff)$, and $(\geff, \DL)$ for simulations based on model
$P_1$ in \cref{fig:comparison-jkdes-corr-exB} and model $P_0$ in
\cref{fig:comparison-jkdes-indep} exhibit qualitative similarities.
Although the qualitative nature of the simulated distributions suggest
no remarkable difference in the physics over the narrow hyperparameter
range, in general the simulations based on model $P_1$ follow a
different sampling procedure than the simulations based on model
$P_0$; in contrast to $P_0$, the causal relationships embedded in
model $P_1$ ensure the pore-scale geometries sampled for the numerical
experiment are always consistent with hierarchical nanoporous material
in \cref{fig:pore-structure} even over extended hyperparameter ranges
such as \cref{tab:extended-range}.

A quantitative comparison suggests that the difference in sampling
procedures leads to distinct distributions with statistical
significance. To compare the densities in
\cref{fig:comparison-kdes-indep,fig:comparison-kdes-corr-exB} and in
\cref{fig:comparison-jkdes-indep,fig:comparison-jkdes-corr-exB}
quantitatively, we work directly with samples from the gPCEs,
employing a two-way statistical test on the equality of distributions.
Since the data in
\cref{fig:comparison-kdes-corr-exB,fig:comparison-kdes-indep} and in
\cref{fig:comparison-jkdes-corr-exB,fig:comparison-jkdes-indep} appear
commensurable, we select a nonparametric Cramér test
(\cite{BaringhausFranz:2004tt}) indicated to be sensitive against
location alternatives that is applicable to both univariate and
multivariate distributions so as to have a consistent presentation.
Although the estimated densities in
\cref{fig:comparison-kdes-corr-exB} and in
(\cref{fig:comparison-jkdes-corr-exB}) are superficially similar to
the densities in \cref{fig:comparison-kdes-indep} (respectively,
\cref{fig:comparison-jkdes-indep}), the statistical tests each based
on \num{2000} sample values, summarized in
\cref{tab:cramer-test-univariate}, reject the hypothesis on equality
of distributions with high statistical significance for all but one
comparison.

\begin{table}[!h]
  \centering
  \caption{Results for a two-way nonparametric Cramér test
    (\cite{BaringhausFranz:2004tt}) on the hypothesis of equality of
    the empirical distirubtions for comparable variables displayed in
    \cref{fig:comparison-kdes-indep,fig:comparison-kdes-corr-exB,fig:comparison-jkdes-indep,fig:comparison-jkdes-corr-exB}
    each based on \num{2000} values sampled using a gPCE.}
  \label{tab:cramer-test-univariate}
  \begin{tabular}{lSSSS[table-format = <1.3]cl}
    \toprule
    Variable & \multicolumn{1}{c}{Cramér-statistic} & \multicolumn{1}{c}{Critical value} 
    & \multicolumn{1}{c}{Conf.\ interval} & \multicolumn{1}{c}{$p$-value} 
    & \multicolumn{1}{c}{Result} & Figures \quad \quad\\ \midrule
    $\DL$& 10.16& 0.3116& 0.95& < 0.001& reject
                                 & \rdelim\}{3}{0mm}[\ref{fig:comparison-kdes-indep} vs.\ \ref{fig:comparison-kdes-corr-exB}]\\
    $\DT$& 1.154& 0.1049& 0.95& < 0.001& reject & \\
    $\geff$& 0.2232& 0.3835& 0.95& 0.154& accept & \\ \midrule

    $(\DL,\DT)$& 10.2& 0.3154& 0.95& < 0.001& reject
                                 & \rdelim\}{3}{0mm}[\ref{fig:comparison-jkdes-indep} vs.\ \ref{fig:comparison-jkdes-corr-exB}]\\
    $(\DT,\geff)$& 1.3& 0.3854& 0.95& < 0.001& reject & \\
    $(\geff,\DL)$& 7.525& 0.4533& 0.95& < 0.001& reject & \\ %  \midrule
 
    \bottomrule
  \end{tabular}
\end{table}

Using model $P_1$, one is able to consider numerical experiments over
the physical range of hyperparameters in \cref{tab:extended-range}
that lead to non-trivial correlations as demonstrated in
\cref{fig:corr-exB-extended}. Over the physical range of
hyperparameters, the marginal and joint densities in
\cref{fig:comparison-marginal-kdes-exBext,fig:comparison-jkdes-corr-exBext},
respectively, for the Darcy-scale flow variable are markedly different
from the marginal and joint densities observed in
\cref{fig:comparison-kdes-indep,fig:comparison-kdes-corr-exB,fig:comparison-jkdes-indep,fig:comparison-jkdes-corr-exB}.
In comparing \cref{fig:comparison-marginal-kdes-exBext} to
\cref{fig:comparison-kdes-corr-exB}, we observe that density for the
effective rate constant $\geff$ becomes more positively skewed and
more peaked suggesting less variance in the estimate of $\geff$ over
the physical hyperparameter range. In contrast, we observe in
comparing \cref{fig:comparison-marginal-kdes-exBext} to
\cref{fig:comparison-kdes-corr-exB} that the density for $\DT$ becomes
more uniform in distribution (and more like the distribution of $\DL$)
suggesting that the simulations over the physical hyperparamter range
realize a richer variety of transverse diffusions. Similarly, in
\cref{fig:comparison-jkdes-corr-exBext} we observe that joint
distributions for the Darcy-scale flow variables employing model $P_1$
over the physical hyperparameter range demonstrate similar
qualitative changes, with the joint densities involving $\geff$
narrowing in variability and realizing a more variety in the observed
transverse diffusion $\DT$. Altogether, the difference between the
densities for the Darcy-scale flow QoIs for model $P_1$ suggest that
different physics is observed over the narrow versus the physical
hyperparameter ranges. In the present context, these differences in
physics can have a profound impact on decision support tasks such as
experimental design thereby underscoring the importance of including
causal relationships in mathematical and statistical modeling
processes. In the section that follows, we continue to investigate the
impact of causality on uncertainty in Darcy flow variables and QoIs by
investigating methodologies for global sensitivity analysis.

\begin{remark}
  In any of the Bayesian network representations
  \cref{eq:P_0,eq:P_1,eq:P_2,eq:dist-DL}, uncertainty and error in the
  homogenization can be included by putting a distribution on $\B{X}$
  that is not trivial, i.e.\ by replacing \cref{eq:dist-X} with a
  distribution that captures error in the homogenization map or that
  compares distributions resulting from different upscaling
  techniques. Uncertainty and error in other relevant processes, such
  as the choice of the KDE, might also be incorporated into the
  Bayesian network for the full statistical model and analyzed
  similarly. Thus, the Bayesian network PDE formulation provides a
  systematic framework for building complete predictive models
  including forward physical models, transitions between scales, and
  uncertainties in parameters, mechanisms and parameter or model
  constraints.
\end{remark}

\section{Global sensitivity analysis and effect ranking}
\label{sec:global-sa}

Recall that we are interested in simulating porous media to inform
general decision tasks concerning the design of materials with
targeted macroscopic properties, such as effective diffusions and
sorption constants, through engineered microscopic pore-scale
structures. In this context, it is important to analyze the
sensitivity of macroscopic QoIs with respect to uncertainties in
pore-scale properties. Simulations of macroscopic material performance
that are highly sensitivity to distributional changes in microscopic
pore-scale processes and structures would undermine the generality of
such investigations.

In general, sensitivity analysis is a key component of uncertainty
quantification and informs decision tasks such optimal experimental
design and the analysis of model robustness, identifiability, and
reliability
(\cite{KomorowskiEtAl:2011sa,SaltelliEtAl:2000sa,SaltelliEtAl:2000sb}).
Local sensitivity analysis is suited to situations where the
(hyper)parameters are known with some confidence and small
perturbations are relevant. In contrast, in our application of
interest the mapping from pore-scale input distributions to
Darcy-scale flow variables is nonlinear, includes several
computational steps, and we have no \emph{a priori} information on the
form of the model for the Darcy-scale variables (cf.\
\cref{sec:unknown-darcy-var-models}). Therefore, the present
application demands global sensitivity analysis methods that explore
the whole space of uncertain input factors
(\cite{SaltelliEtAl:2008gs}).

Variance-based sensitivity analysis methods, such as Sobol' indices
\cite{Sobol:1993sa} and total sensitivity indices
\cite{HommaSaltelli:1996im}, rank input factors and higher order
interactions of input factors in terms of their contributions to the
variance of a QoI. In particular, Sobol' sensitivity indices are used
in \cite{UmZhangKatsoulakisEtAl:2017aa} to analyze the global
sensitivity of Darcy-scale QoIs to first and second order interactions
among input parameters for the multiscale porous media model
considered here under the assumption of independent uninformed priors
on the pore-scale distributions. However, such variance-based methods
for assessing global sensitivity are not easily applied or interpreted
in the case of the dependent input parameters introduced through
causal relationships outlined in \cref{sec:prob-graph-model}.
Moreover, we observe that the distributions for macroscopic
Darcy-scale flow variables exhibit non-Gaussian behavior (cf.\
marginal and joint densities in
\cref{fig:comparison-kdes,fig:comparison-joint-kdes} in
\cref{sec:quant-uncert-darcy}). Methods that rely on moment
information alone may be insufficient to capture the full complexity
of these interactions.

\subsection{Mutual information for global sensitivity analysis}

Presently, we employ global sensitivity indices based on information
theoretic concepts (\cite{CoverThomas:2006in,Soofi:1994if}) that rely
on empirical distributions as opposed to moments. There is a rich
literature on moment-independent indices for local sensitivity
analysis
(\cite{PantazisEtAl:2013rn,PantazisKatsoulakis:2013re,KomorowskiEtAl:2011sa,MajdaGershgorin:2011sa,MajdaGershgorin:2010qu})
as well as global sensitivity analysis
(\cite{ChunEtAl:2000im,Borgonovo:2006uq,LiuEtAl:2006re,Borgonovo:2007sa,LiuHomma:2009uq,CastaingsEtAl:2012sa,Rahman:2016fs}).
In particular, global sensitivity indices based on various information
theoretic notions are well established in the literature
(\cite{ParkEtAl:1994sa,LiuEtAl:2006re,LudtkeEtAl:2008sa,VetterTaflanidis:2012sa}).
In \cite{LudtkeEtAl:2008sa}, discrete mutual information based
sensitivity indices are demonstrated as an effective tool for discrete
probability distributions arising in biochemical reaction networks in
systems biology. For our application of interest, sensitivity indices
and rankings based on the \emph{differential mutual information}
provide a suitable measure of effect that overcomes the twin
challenges of causally related inputs and non-Gaussian QoI. The
differential mutual information has explicit connections to more
general information theoretic concepts and we provide interpretations
of these in the context of uncertainty quantification and global
sensitivity analysis.

The differential mutual information between continuous random
variables $V$ and $W$,
\begin{equation}
  \label{eq:mutual-info}
  I(V; W) \defeq \iint \log\left(\frac{f_{V,W}(v,w)}{f_V(v)f_W(w)} \right) 
  f_{V,W}(v,w) \dd v \dd w \,,
\end{equation}
quantifies the statistical dependence, that is, the amount of shared
information, between $V$ and $W$ provided that all of the densities
above exist and the marginals are non-zero. Importantly,
\cref{eq:mutual-info} applies to dependent random variables $V$ and
$W$, such as random variables that share a causal relationship, and we
observe that $I(V,W) = 0$ if $V$ and $W$ are independent, i.e., if
$\Pb(V,W) = \Pb(V)\Pb(W)$. Further, \cref{eq:mutual-info} applies to
random variables with very general marginal and joint distributions
including distributions that are non-Guassian. The differential mutual
information is precisely the relative entropy $\RE$ (or
Kullback--Leibler divergence) between the joint and marginal
densities,
\begin{equation}
  \label{eq:mutual-info-KL-div}
  I(V; W) = \RE ( \Pb(V,W) \midbars \Pb(V) \Pb(W) )\,,
\end{equation}
a well known pseudo-distance used in variational inference
(\cite{WainwrightJordan:2008ef,BleiKucukelbirMcAuliffe:2017vi}),
machine learning (\cite{Bishop:2006ml}), and model selection
(\cite{BurnhamAnderson:2002ms}) as well as other areas. Lastly, the
differential mutual information $I(V,W)$ is the limiting value of the
discrete mutual information (the supremum over all partitions of $V$
and $W$) and therefore shares all of the same properties as its
discrete counterpart (\cite{CoverThomas:2006in}); in particular, we
will compute the differential mutual information using empirical
distributions.

\subsection{Estimators for mutual information global sensitivity
  indices}

For each QoI $g$, we consider a global sensitivity index,
\begin{equation}
  \label{eq:mutual-info-sensitivity-index}
  S_{\Theta}(g) \defeq I(g; \Theta)\,,
\end{equation}
based on the mutual information between a Darcy-scale QoI $g$ and each
uncertain pore-sclae input parameter
$\Theta \in \B{\Theta} = (\Theta_1, \dots, \Theta_{N_p})$.
Intuitively, the index \cref{eq:mutual-info-sensitivity-index}
measures the predictability of $g$ given knowledge of $\Theta$ through
the discrepancy between the joint density and the product of the
marginal densities appearing in \cref{eq:mutual-info}.

We estimate the sensitivity index
\cref{eq:mutual-info-sensitivity-index} using a MC approximation that
relies on empirically estimated distributions, i.e.\ KDEs for the
density functions of Darcy-scale QoIs. Indeed, a benefit to using the
mutual information is the availability of methods relying on plug-in
estimators for \cref{eq:mutual-info} with corresponding numerical
analysis (including
\cite{KraskovEtAl:2004mi,Paninski:2003mi,AntosKontoyiannis:2001cp}).
Recall that we obtain approximate densitities $\bar{f}$ using
univariate and multivariate KDEs, in \cref{eq:general-kde-g} and
\cref{eq:general-jkde-g1-g2}, respectively, that are in turn obtained
using the gPCE surrogates $\hat{g}$ defined in \cref{eq:gPCe-g}. Due
to the availability of these surrogates, we do not sample input-output
pairs as suggested by the form of \cref{eq:mutual-info} but instead
consider the equivalent representation,
\begin{equation}
  \label{eq:mutual-info-marginal-sampling}
  I(V; W) = \iint \log\left(\frac{f_{V,W}(v,w)}{f_V(v)f_W(w)} \right) 
  \frac{f_{V,W}(v,w)}{f_V(v)f_W(w)} f_{V}(v) f_{W}(w) \dd v \dd w \,.
\end{equation}
Thus we compute the statistical estimator
$\widehat{S}_\Theta(g) \approx S_\Theta(g)$,
\begin{equation}
  \label{eq:mc-estimator-S}
  \widehat{S}_\Theta (g) \defeq 
  \frac{1}{N} \sum_{k=1}^{N} 
  \log \left(\frac{\bar{f}_{\hat{g},\Theta} (\hat{g}^k, \Theta^k)}{\bar{f}_{\hat{g}} (\hat{g}^k) f_{\Theta}(\Theta^k)} \right) \frac{\bar{f}_{\hat{g},\Theta} (\hat{g}^k, \Theta^k)}{\bar{f}_{\hat{g}} (\hat{g}^k) f_{\Theta}(\Theta^k)} \,,
\end{equation}
based on a MC approximation of \cref{eq:mutual-info-marginal-sampling}
using KDEs as plug-in estimates where appropriate.

\begin{figure}[!h]
  \centering \includegraphics[width=0.4\textwidth]{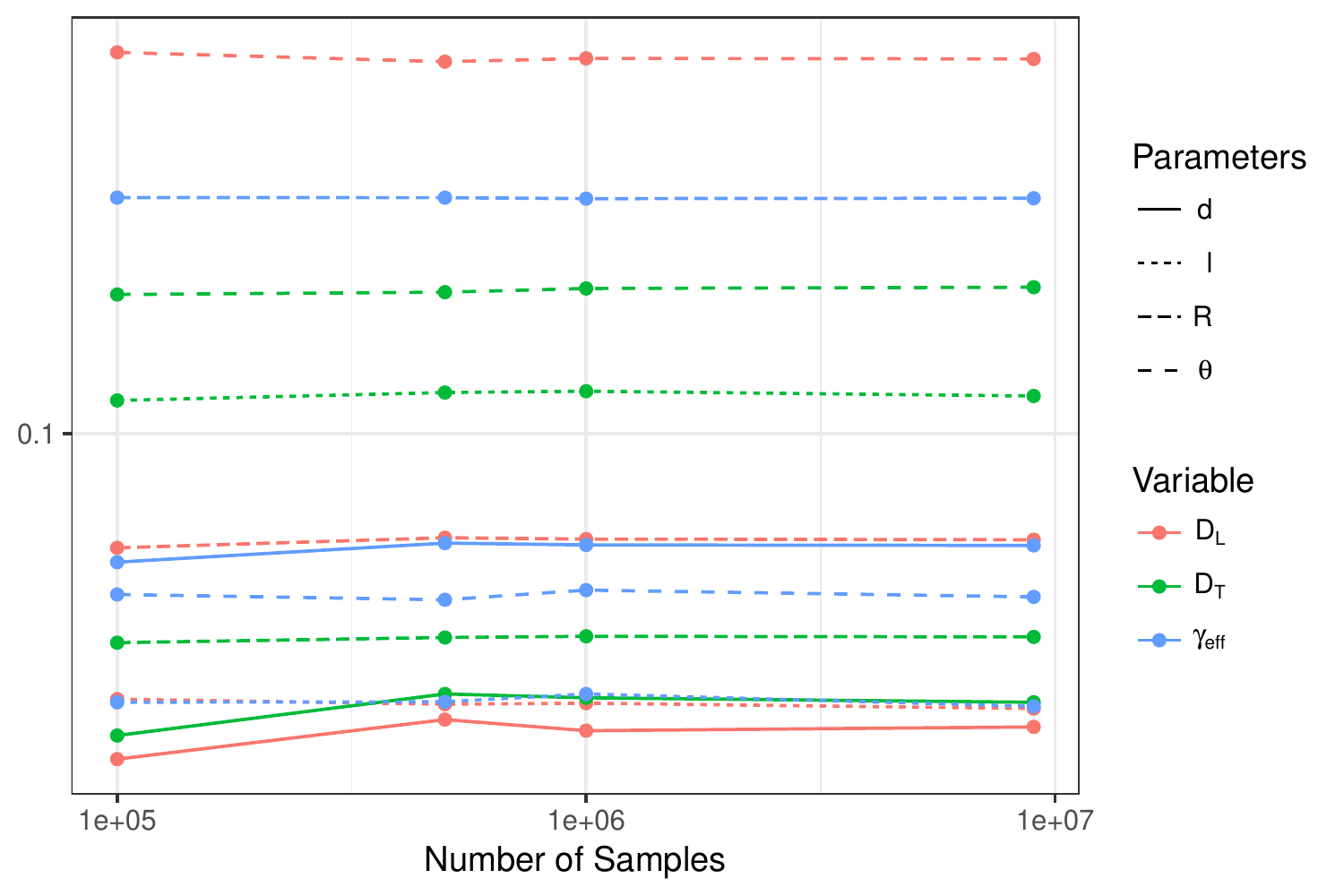}
  \caption{The convergence of the MC estimator
    \cref{eq:mc-estimator-S} for model $P_1$ in \cref{eq:P_1} over the
    narrow hyperparameter range in \cref{tab:comparison-range} for
    each Darcy-scale flow variable for each input parameter
    demonstrates that \num{1e5} samples is sufficient for the
    numerical experiments.}
  \label{fig:mc-convergence}
\end{figure}

In the left-hand side of \cref{tab:mutual-info-and-ranking}, we report
the value of the statistical estimator $\widehat{S}_\Theta (g)$
related to the numerical experiments presented in
\cref{sec:Darcy-scale-qoi} concerning the probabilistic model $P_1$ in
\cref{eq:P_1} with causal inputs (see \cref{fig:conditional-d-and-l})
over both the narrow and physical hyperparameter ranges in
\cref{tab:comparison-range,tab:extended-range}. We are interested in
the global sensitivity of the Darcy-scale QoIs $g = \DL$, $g=\DT$, and
$g=\geff$ with respect to each of input parameters
$\Theta \in \B{\Theta} = (\Theta_R, \Theta_\theta, \Theta_d,
\Theta_l)$ representing the pore-scale features identified in the
hierarchical nanoporous material in \cref{fig:pore-structure}. The
computed estimators \cref{eq:mc-estimator-S} are based on first
generating random variables $\Theta^k$ and $\hat{g}^k$ for
$k = 1,\dots, N = \num{1e7}$ using the respective gPCE surrogate
obtained by the workflow outlined in \cref{sec:quant-uncert-darcy}.
These samples are used to form the respective KDE on an $\eta$-grid of
size 128 by 128 which in turn are then used to form the plug-in
quantity
\begin{equation*}
  X^j = \log\left(\frac{\bar{f}_{\Theta,\hat{g}}(\Theta^j, \hat{g}^j)}
    {\bar{f}_\Theta(\Theta^j) \bar{f}_{\hat{g}}(\hat{g}^j)}\right)
  \frac{\bar{f}_{\Theta,\hat{g}}(\Theta^j, \hat{g}^k)}
  {\bar{f}_\Theta(\Theta^j) \bar{f}_{\hat{g}}(\hat{g}^j)}\,,
\end{equation*}
where we define $\log(\frac{0}{0}) \cdot \frac{0}{0} = 0$ (re-using
$\num{1e5}$ of the samples $(\Theta^j, \hat{g}^j)$ generated during
density estimation). The sensitivity index is then estimated using the
MC estimator,
$\widehat{S}_\Theta(\hat{g}) = \frac{1}{M} \sum_{j=1}^{M} X^j$ with
$M = \num{1e5}$ samples. We observe that the $M = \num{1e5}$ samples
utilized for the MC estimator is suggested to be sufficient by the
convergence demonstrated in \cref{fig:mc-convergence} for the
experiment with model $P_1$ over the range
\cref{tab:comparison-range}; similar convergence observations were
made for the other experiments.

\begin{table}[!h] \centering
  \caption{The global sensitivity index $S$ in
    \cref{eq:mutual-info-sensitivity-index} based on the mutual
    information quantifies the effect of pore-scale uncertainties in
    terms of how additional knowledge of the input
    $\Theta \in \B{\Theta}$ reduces uncertainty in our prediction of
    the macroscopic Darcy-scale variables $\DL$, $\DT$, and $\geff$.
    Below, the estimator $\widehat{S}$ in \cref{eq:mc-estimator-S} and
    the associated ranking $\widehat{r}$ in
    \cref{eq:mutual-info-ranking} is given for the probabilistic
    models $P_1$ in \cref{eq:P_1} over both the hyperparameter ranges
    in \cref{tab:comparison-range} and \cref{tab:extended-range} (cf.\
    graphical representation of rankings in
    \cref{fig:mutual-info-ranking}).}
  \label{tab:mutual-info-and-ranking}
  \begin{tabular}{@{}llllllllll@{}} \toprule &&
    \multicolumn{3}{c}{Mutual information \cref{eq:mc-estimator-S}} &&
    \multicolumn{3}{c}{Ranking effects \cref{eq:mutual-info-ranking}}
    & Hyperparameters\\
    $\Theta$ && \multicolumn{1}{c}{$\widehat{S}_\Theta(\DL)$} &
                                                                \multicolumn{1}{c}{$\widehat{S}_\Theta(\DT)$}
    & \multicolumn{1}{c}{$\widehat{S}_\Theta(\geff)$} & &
                                                          \multicolumn{1}{c}{$\widehat{r}_\Theta(\DL)$}
    & \multicolumn{1}{c}{$\widehat{r}_\Theta(\DT)$}
                                                      & \multicolumn{1}{c}{$\widehat{r}_\Theta(\geff)$}\\
    \midrule
    $\Theta_R$ & %$P_1$, \cref{tab:comparison-range}
                                                                    & \num{0.0419}& \num{0.0424}& \num{0.8955}& 
    & \num{0.0724}& \num{0.1177}& \num{0.8509}& 
                                                \rdelim\}{4}{0mm}[\cref{tab:comparison-range}]\\
    $\Theta_\theta$ & %$P_1$, \cref{tab:comparison-range}
                                                                    & \num{0.5074}& \num{0.2049}& \num{0.0366}& 
    & \num{0.8770}& \num{0.5692}& \num{0.0348}& \\
    $\Theta_d$ & %$P_1$, \cref{tab:comparison-range}
                                                                    & \num{0.0150}& \num{0.0878}& \num{0.0271}& 
    & \num{0.0259}& \num{0.2438}& \num{0.0258}& \\
    $\Theta_l$ & %$P_1$, \cref{tab:comparison-range}
                                                                    & \num{0.0143}& \num{0.0249}& \num{0.0932}& 
    & \num{0.0247}& \num{0.0692}& \num{0.0885}& \\

    \midrule
    $\Theta_R$ & %$P_1$, \cref{tab:extended-range}
                                                                    & \num{0.0655}& \num{0.0714}& \num{0.2878}& 
    & \num{0.1354}& \num{0.2470}& \num{0.5242}& 
                                                \rdelim\}{4}{0mm}[\cref{tab:extended-range}]\\
    $\Theta_\theta$ & %$P_1$, \cref{tab:extended-range}
                                                                    & \num{0.3539}& \num{0.0312}& \num{0.0251}& 
    & \num{0.7312}& \num{0.1079}& \num{0.0457}& \\
    $\Theta_d$ & %$P_1$, \cref{tab:extended-range}
                                                                    & \num{0.0225}& \num{0.1646}& \num{0.0823}& 
    & \num{0.0465}& \num{0.5697}& \num{0.1499}& \\
    $\Theta_l$ & %$P_1$, \cref{tab:extended-range}
                                                                    & \num{0.0420}& \num{0.0218}& \num{0.1539}& 
    & \num{0.0868}& \num{0.0754}& \num{0.2802}& \\

    \bottomrule
  \end{tabular}
\end{table}

\begin{remark}
  The differential mutual information \cref{eq:mutual-info} can be
  expressed as an expected value of the relative entropy between
  conditional and marginal distributions,
  \begin{equation}
    \label{eq:mutual-info-input-output-pair}
    I(g; \Theta) 
    = \E_{\Theta} \left[ \RE (P(g\given \Theta) \midbars P(g))\right]\,,
  \end{equation}
  for the input-output pair $(g(Y(\B{\Theta})),\Theta)$ provided
  $\Pb(g, \Theta) = \Pb(g, \given \Theta) \Pb(\Theta)$, i.e., all the
  densities exist. In \cite{Rahman:2016fs}, a family of sensitivity
  measures is given by replacing $\RE$ in
  \cref{eq:mutual-info-input-output-pair} with a general class of
  Csiszár $\varphi$-divergences
  (\cite{AliSilvey:1966fd,Csiszar:1967fd,LieseVajda:2006it}); here we
  demonstrate an implementation using gPCE surrogates that requires
  the representation \cref{eq:mutual-info-marginal-sampling} and
  restrict our attention to the differential mutual information owing
  to the clear interpretation as a measure of effect in terms of
  statistical dependence and shared information.
\end{remark}

\begin{remark}
  One can also consider higher order effects that include interactions
  between a subset of parameters using the conditional differential
  mutual information that takes the form of conditional expectations
  of the relative entropy between joints and marginals, as in
  \cite{LudtkeEtAl:2008sa}. These higher order interactions can be
  interpreted similarly in terms of dependence and shared information.
\end{remark}

\subsection{Ranking impact of uncertainty in correlated pore-scale
  inputs on Darcy-scale QoIs}
\label{sec:ranking-pore-scale}

\begin{figure}[!h]
  \centering
  \begin{subfigure}[t]{0.48\textwidth}
    \includegraphics[width=1\linewidth]{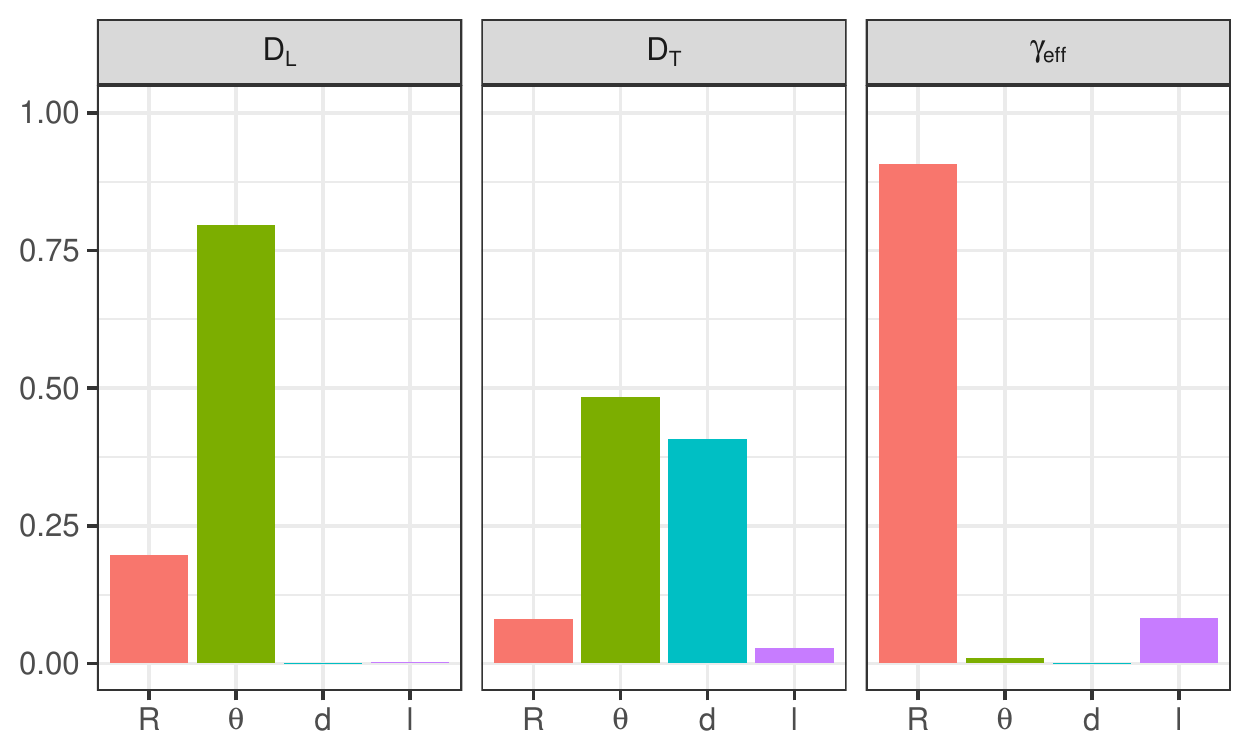}
    \caption{First order Sobol' indices for independent priors model
      $P_0$ over narrow hyperparameter range (reproduced from data in
      \cite{UmZhangKatsoulakisEtAl:2017aa} for the convenience of the
      reader).}
    \label{fig:sobol-ranking}
  \end{subfigure} \quad
  \begin{subfigure}[t]{0.48\textwidth}
    \includegraphics[width=1\linewidth]{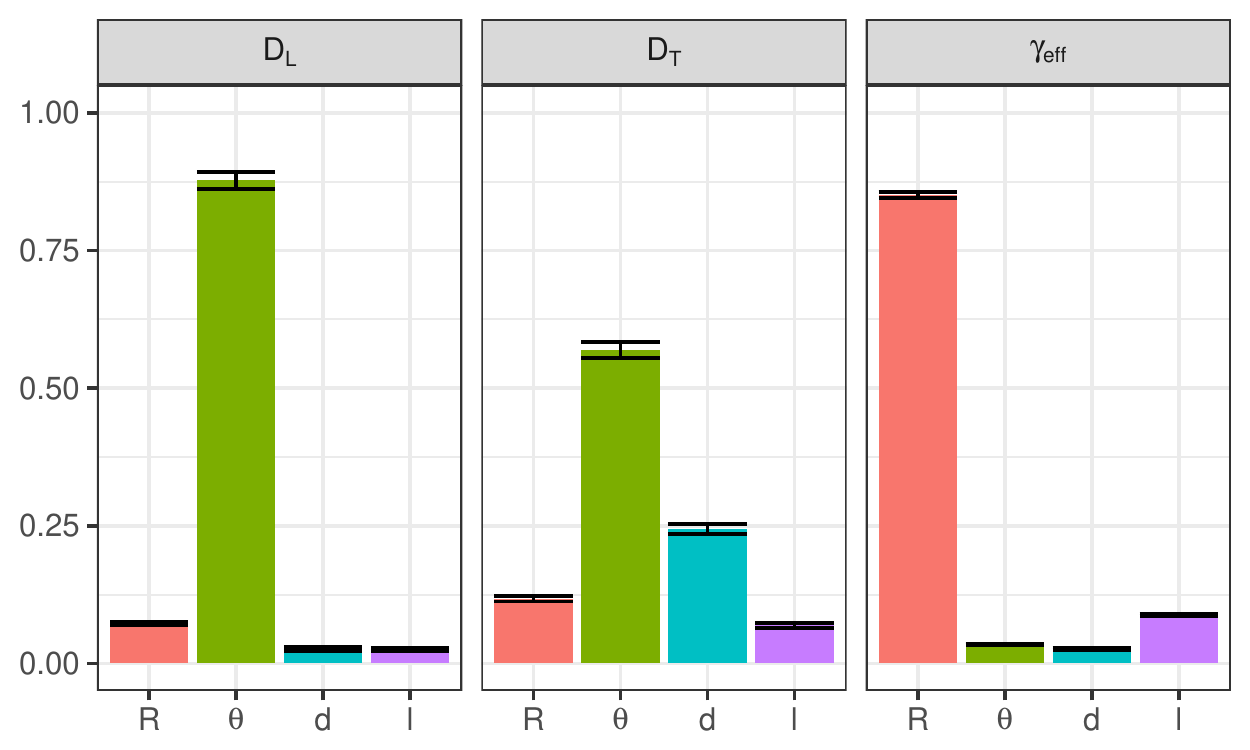}
    \caption{Ranking effects $\widehat{r}_\Theta (\hat{g})$ for causal
      model $P_1$ over narrow hyperparameter range.}
    \label{fig:rank-P1-narrow}
  \end{subfigure}\\
  \vspace*{2ex}
  \begin{subfigure}[b]{0.48\textwidth}
    \includegraphics[width=1\linewidth]{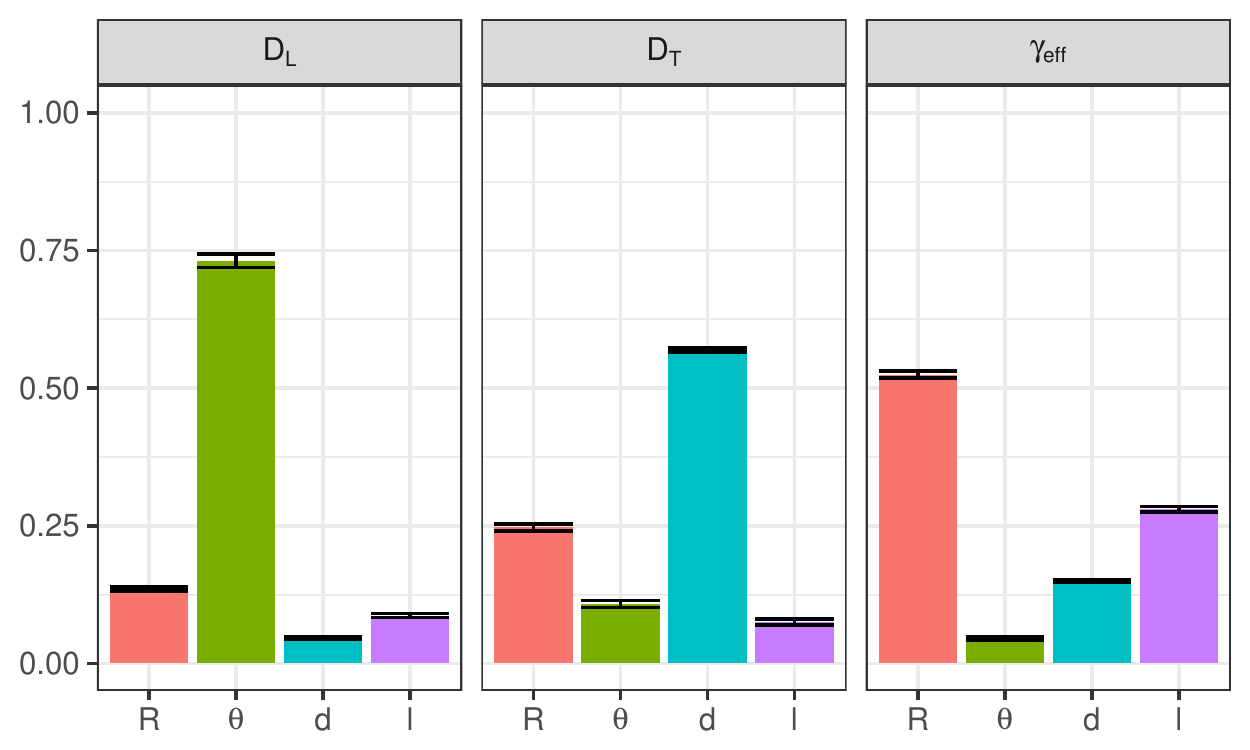}
    \caption{Ranking effects $\widehat{r}_\Theta (\hat{g})$ for causal
      model $P_1$ over physical hyperparameter range.}
    \label{fig:rank-P1-extended}
  \end{subfigure}
  \caption{The rankings in \cref{fig:rank-P1-extended} associated with
    the causal model $P_1$ in \cref{eq:P_1} over the physical
    hyperparameter range in \cref{tab:extended-range} yield parameter
    rankings that are consistent with our understanding of the physics
    of the hierarchical nanoporous material in
    \cref{fig:pore-structure}. In contrast, the Sobol' index rankings
    in \cref{fig:sobol-ranking} for the independent priors model $P_0$
    in \cref{eq:P_0} over the narrow hyperparameter range in
    \cref{tab:comparison-range} suggest that the transverse diffusion
    $\DT$ is most sensitive to the parameter $\theta$, related to the
    angle of overlap between adjacent mesopores in a nanotunnel. The
    mutual information rankings in \cref{fig:rank-P1-narrow} for model
    $P_0$ are consistent with the Sobol' index rankings in
    \cref{fig:sobol-ranking} for model $P_0$ owing to the similarity
    of the correlation structures for pore-scale features over the
    narrow hyperparameter range (cf.\ \cref{fig:correlations}). In
    \cref{fig:rank-P1-narrow,fig:rank-P1-extended}, error bars, that
    indicate the relative error associated with plus/minus two
    standard deviations of the computed sensitivity index, provide an
    indication of confidence in the ranking value.}
  \label{fig:mutual-info-ranking}
\end{figure}

Using the global sensitivity index
\cref{eq:mutual-info-sensitivity-index} we form the ranking,
\begin{equation}
  \label{eq:mutual-info-ranking}
  r_\Theta(g) = \frac{S_\Theta(g)}{\sum_{V \in \B{\Theta}} S_V(g)}\,,
\end{equation}
of the relative contribution of each pore-scale parameter
$\Theta \in \B{\Theta}$ to the global sensitivity of each Darcy-scale
variable $g$. We then obtain the estimate $\widehat{r}_\Theta(g)$
reported in the right-hand side of \cref{tab:mutual-info-and-ranking},
by using the estimator $\widehat{S}$ in \cref{eq:mc-estimator-S} in
place of $S$. The values $\widehat{r}$ are also displayed graphically
in \cref{fig:mutual-info-ranking} with error bars that indicate the
relative error associated with plus or minus two standard deviations
of the computed sensitivity index to provide an indication of
confidence in the ranking.

In \cref{fig:rank-P1-narrow}, we observe that the rankings suggested
by the model $P_1$ in \cref{eq:P_1} with causal inputs in
\cref{fig:conditional-d-and-l} over the narrow hyperparameter range in
\cref{tab:comparison-range} are consistent with the Sobol' index
rankings in \cref{fig:sobol-ranking} that were obtained in
\cite{UmZhangKatsoulakisEtAl:2017aa} for the model $P_0$ with
independent priors. As observed in \cref{sec:Darcy-scale-qoi} this is
to be expected due to the similarity in the correlation structures in
\cref{fig:corr-exInd-narrow,fig:corr-exB-narrow} over the narrow
hyperparameter range. The $\theta$, which is related to the mesopore
radius, is the most influential parameter for both the longitudinal
diffusion $\DT$ and the transverse diffusion $\DT$. The mesopore
radius $R$, which is related to the size of the fluid-solid interface
in \cref{fig:pore-structure}, is the most influential parameter for
the sorption rate constant $\geff$. In general, it is important to
interpret these sensitivity rankings in the context of the
hyperparameter range; the ranking $\widehat{r}_{\Theta_d}$ in
\cref{fig:rank-P1-narrow} is likely to be small as the range for $d$
in \cref{tab:comparison-range} is very narrow (i.e.\ the Darcy-scale
QoIs are insensitive over the narrow range of admissible $d$ values).

In contrast, the rankings in \cref{fig:rank-P1-extended} for the model
$P_1$ over the physical hyperparameter range in
\cref{tab:extended-range} demonstrates considerably different rankings
to \cref{fig:rank-P1-narrow,fig:sobol-ranking} thereby highlighting
once again the impact of causal relationships and the Bayesian network
PDE modeling approach. The rankings in \cref{fig:rank-P1-extended}
indicate that while the $\theta$ and the $R$ are still the most
influential parameters for, respectively, the longitudinal diffusion,
$\DL$, and sorption rate constant, $\geff$, it is $d$, which is
related to the diameter of the nanotube, that is the most influential
parameter for the transverse diffusion $\DT$ through nanotubes. Thus,
for the simple hierarchical nanoporous material in
\cref{fig:pore-structure}, the rankings in \cref{fig:rank-P1-extended}
using model $P_1$ over the physical hyperparameter range reflect our
expectations of the physics better than the experimental observations
with respect to models over the narrow range. Importantly, the
rankings in \cref{fig:mutual-info-ranking} are moment-independent and
therefore suitable for the non-Gaussian behavior of the Darcy-scale
QoIs observed in \cref{fig:comparison-kdes,fig:comparison-joint-kdes}.
Moreover, as the rankings are based on the mutual information, they
can be interpreted as ranking the impact of the pore-scale parameter
(whether correlated or not) on the total uncertainty in the
Darcy-scale quantity of interest.

\section{Alternative probabilistic models reflect different
  engineering and design causal relationships}
\label{sec:alternative-prob-models}

\begin{figure}[!h]
  \centering
  \begin{subfigure}[t]{0.3\textwidth}
    \centering \includegraphics[width=.87\linewidth]{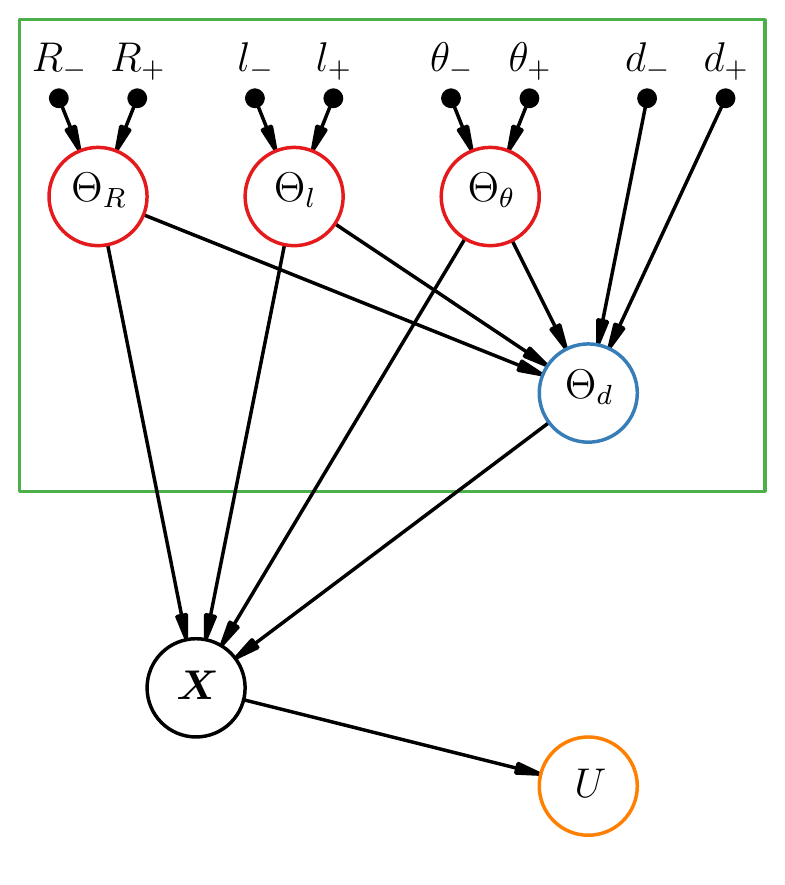}
    \caption{Bayesian network for model $P_2$ encoding different
      causal relationships from $P_1$ (cf.\
      \cref{fig:conditional-d-and-l,fig:general-model-indpriors}).}
    \label{fig:conditional-d-only}%
  \end{subfigure}
  \qquad
  \begin{subfigure}[t]{0.3\textwidth}
    \centering \includegraphics[width=1\linewidth]{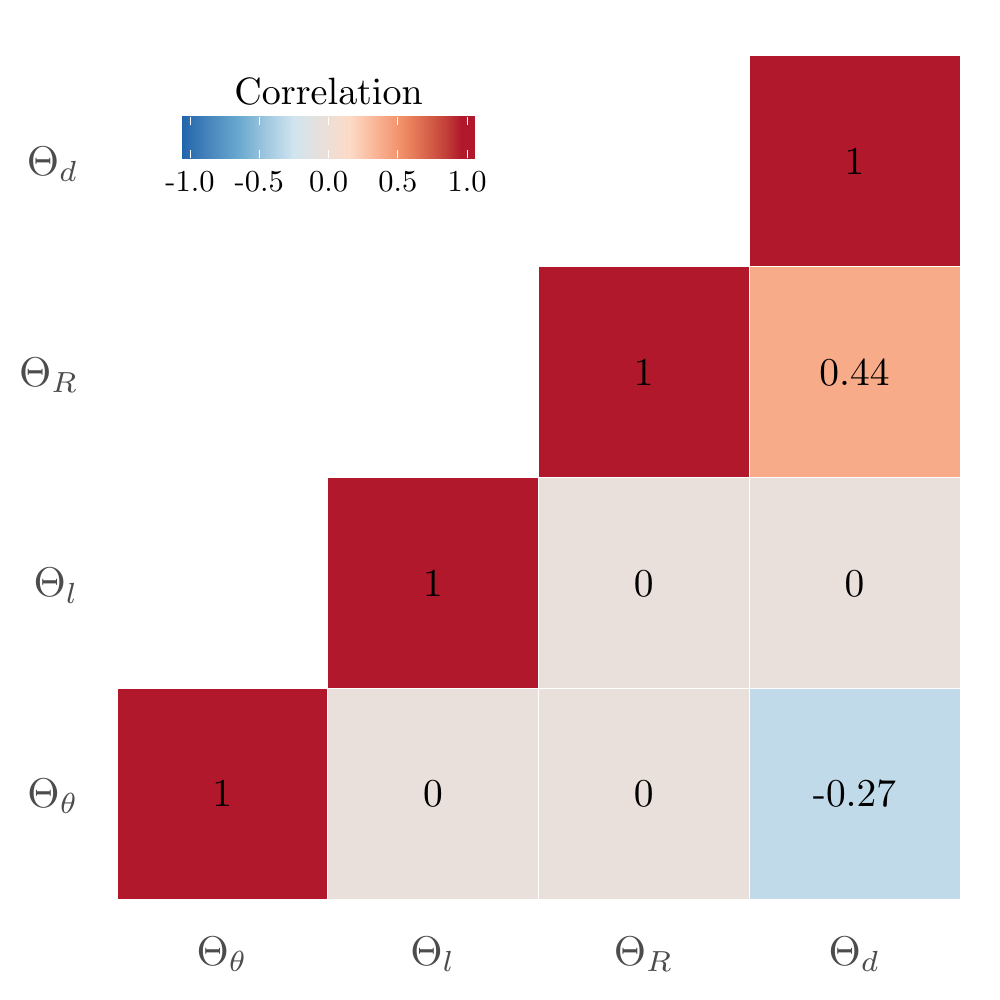}
    \caption{Correlation structure for $P_2$ over physical
      hyperparameter range in \cref{tab:extended-range} (cf.\
      \cref{fig:correlations}).}
    \label{fig:corr-exA-extended}%
  \end{subfigure}
  \qquad
  \begin{subfigure}[t]{0.3\textwidth}
    \centering \includegraphics[width=1\linewidth]{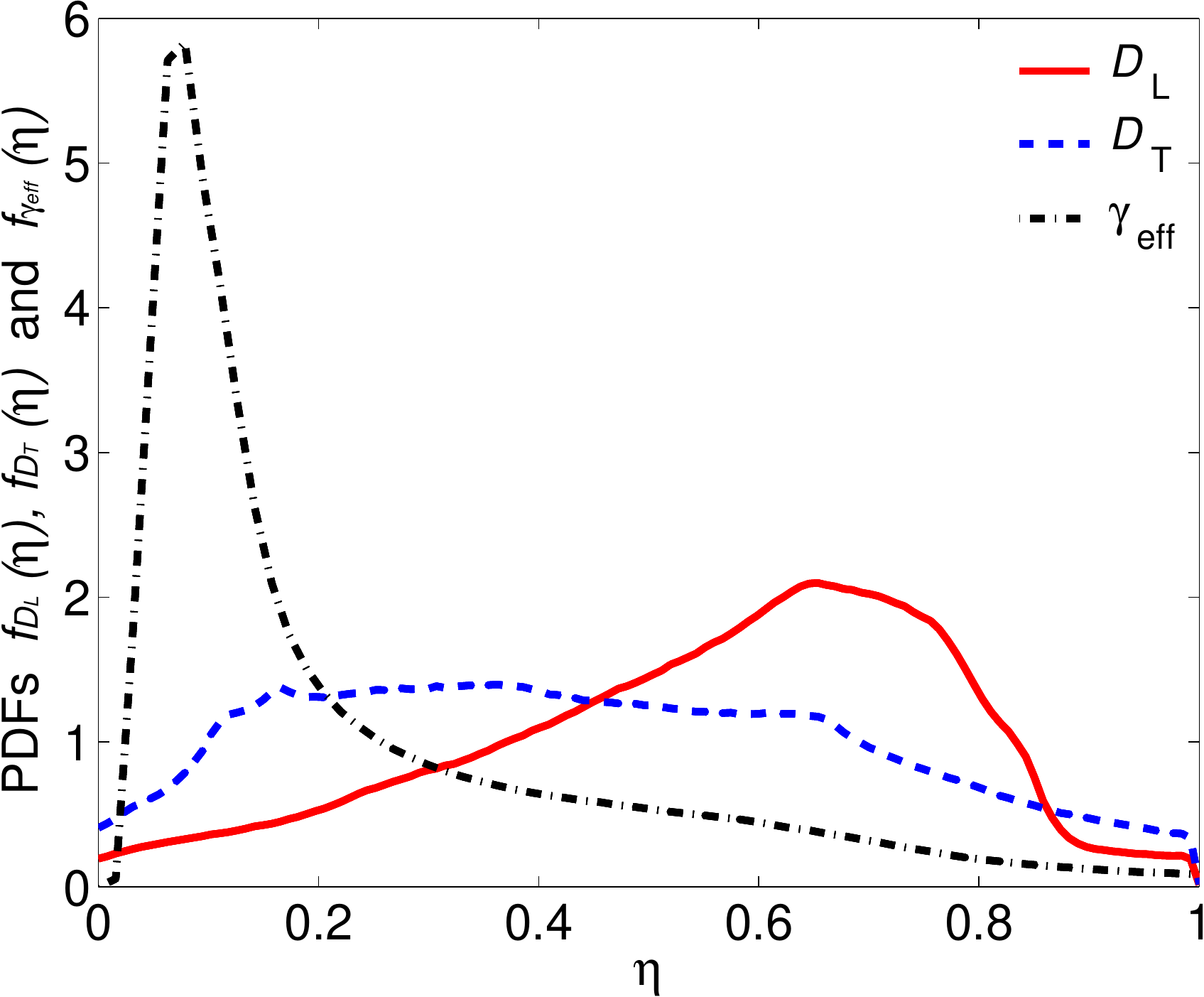}
    \caption{Marginal KDEs for Darcy-scale flow variables for model
      $P_2$ over the physical hyperparameter range (cf.\
      \cref{fig:comparison-kdes}).}
    \label{fig:comparison-marginal-kdes-exAext}%
  \end{subfigure}\\
  \vspace*{2ex}
  \begin{subfigure}[b]{1\textwidth}
    \includegraphics[width=.33\linewidth]{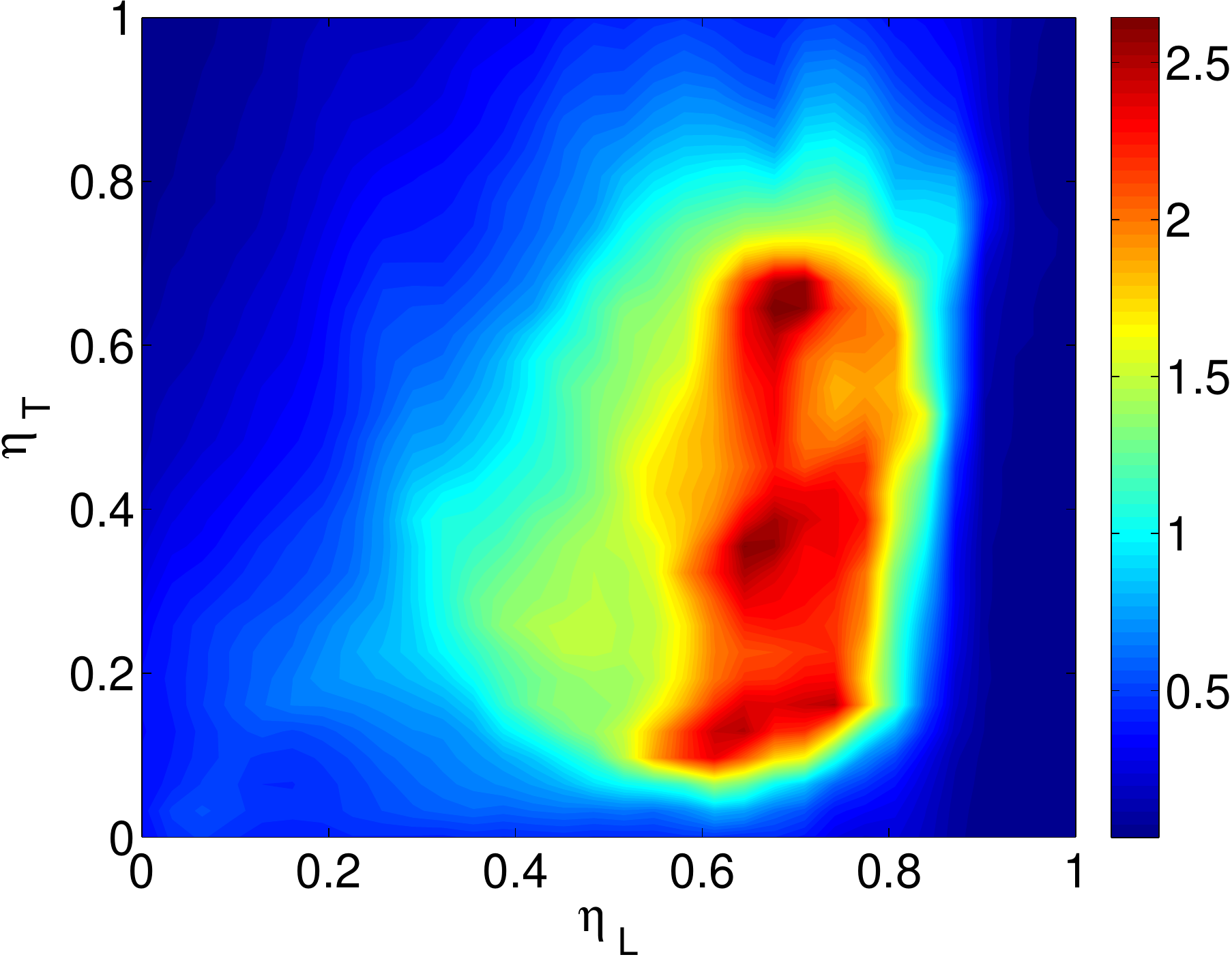}
    \includegraphics[width=.33\linewidth]{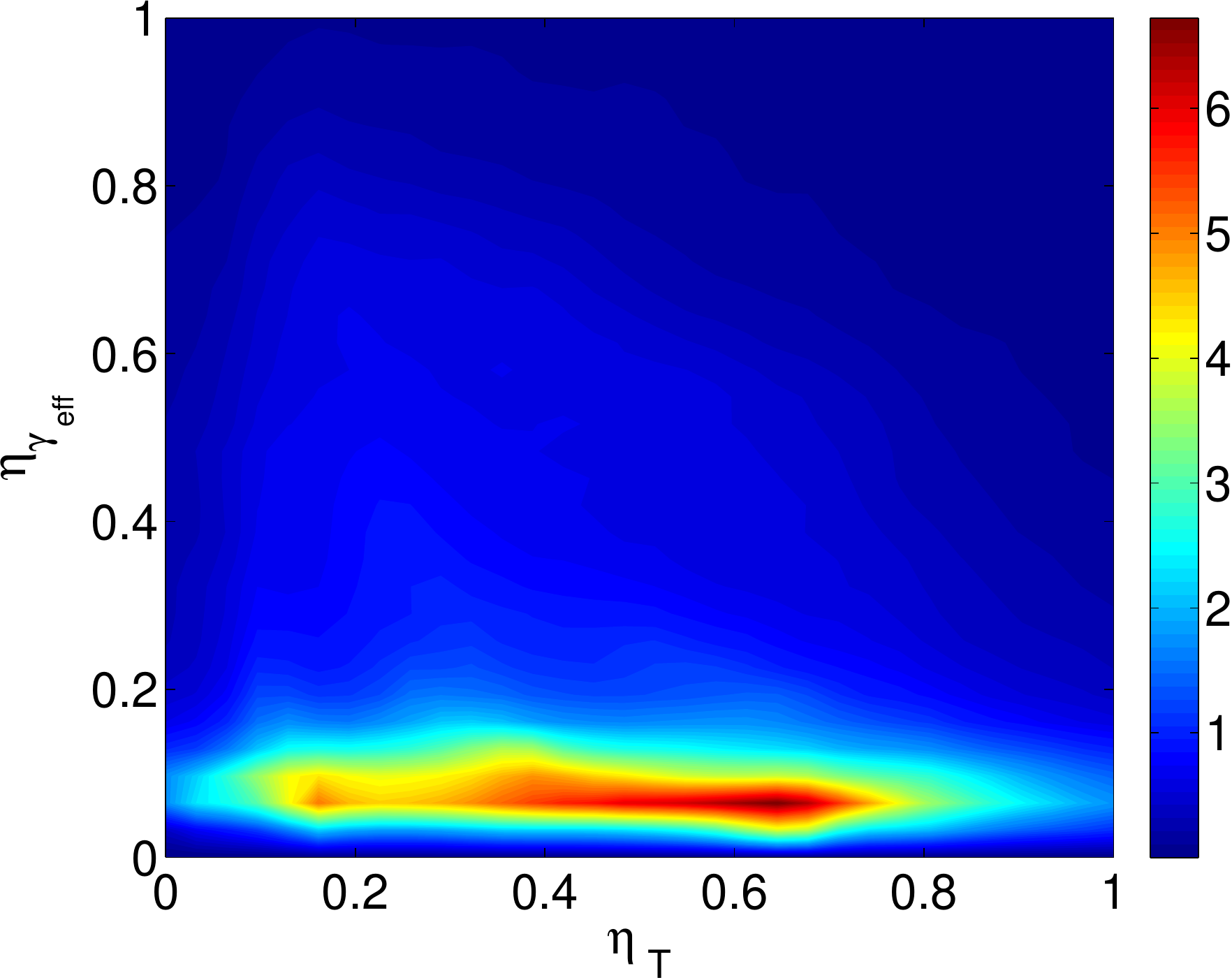}
    \includegraphics[width=.33\linewidth]{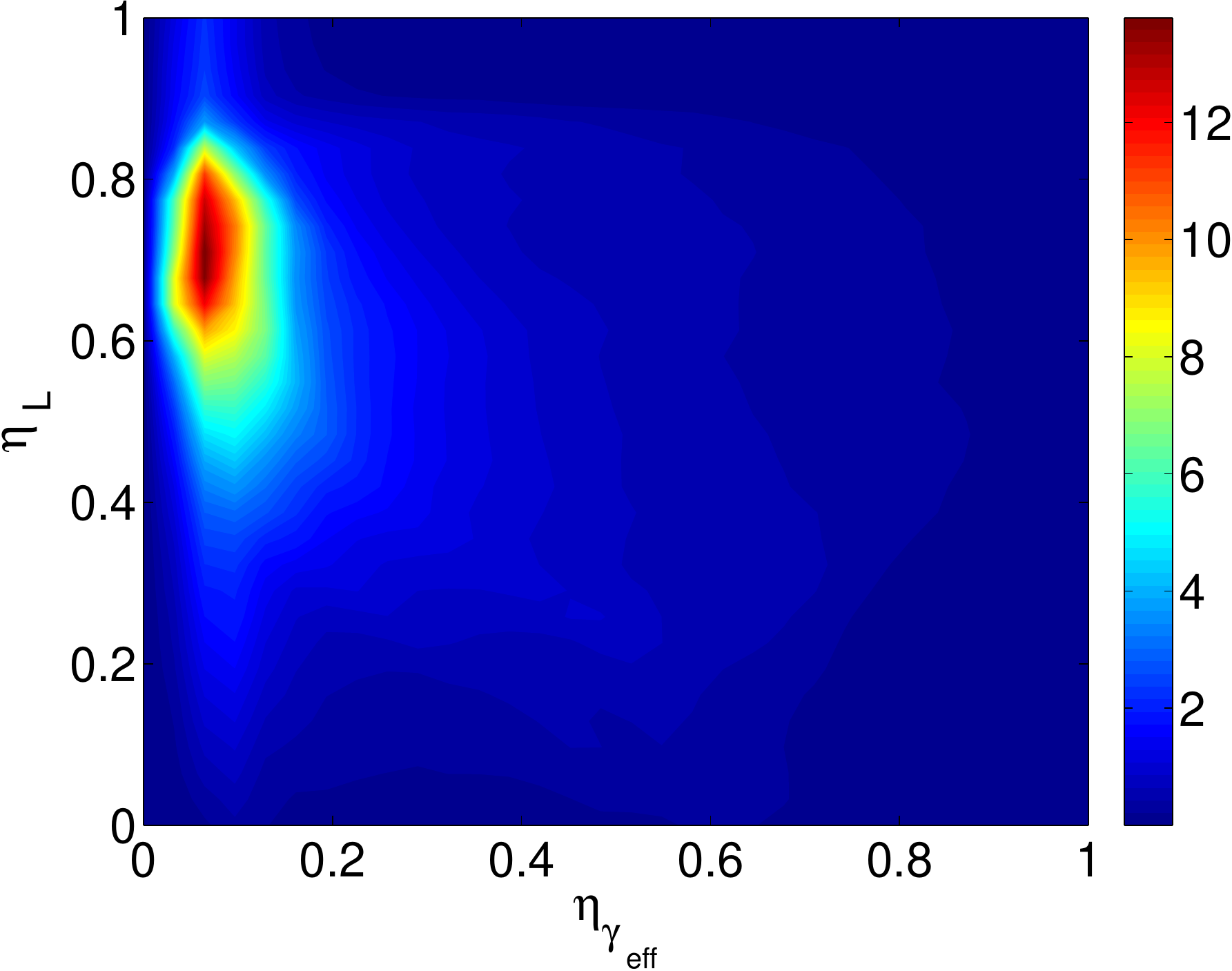}
    \caption{Causal priors (model $P_2$) over the physical
      hyperparameter range (cf.\ \cref{fig:comparison-joint-kdes}).}
    \label{fig:comparison-jkdes-corr-exAext}
  \end{subfigure}
  \caption{An alternative model $P_2$ in \cref{eq:P_2} encodes causal
    relationships that are different from model $P_1$ in \cref{eq:P_1}
    but also respects structural constraints realized by the
    hierarchical nanoporous material in \cref{fig:pore-structure}. The
    model $P_2$ is viewed as mirroring engineering processes and
    design work-flows that are distinct from those of $P_1$.}
  \label{fig:P2}
\end{figure}

We view the construction of the full statistical model for the
multiscale porous media system as mirroring engineering processes and
design work-flows. From this perspective, it may be desirable to build
models containing different causal relationships among the pore-scale
parameters than those previously considered. Recall that the model
$P_1$ in \cref{eq:P_1} with causal priors given by the Bayesian
network in \cref{fig:conditional-d-and-l} is related to the design of
nanotunnels/mesopores, i.e.\ the pore-scale features $R$ and $\theta$
in the hierarchical nanoporous media in \cref{fig:pore-structure}. In
contrast, if it is desirable or possible to choose three aspects, such
as the features $R$, $\theta$, and $l$, independently while
constraining the only remaining parameter, $d$, then one obtains a
second model represented by the Bayesian network in
\cref{fig:conditional-d-only}.

The Bayesian network in \cref{fig:conditional-d-only} corresponds to
placing all of the constraints on $\Theta_d$. That is, choosing a
subset of independent parameters
$\{\Theta_R, \Theta_\theta, \Theta_l\}$ and assuming uniform priors,
\begin{equation}
  \label{eq:Theta_R-Theta_t-Theta_l}
  \Theta_R 
  \sim \unif(R_{-}, R_{+})
  \quad \text{and} \quad 
  \Theta_\theta
  \sim \unif(\theta_{-}, \theta_{+})
  \quad \text{and} \quad 
  \Theta_l \sim \unif(l_{-}, l_{+})\,,
\end{equation}
constrains the distribution of $\Theta_d$. Specifically, in order for
the choice of $d$ to be consistent with the features depicted in
\cref{fig:pore-structure}, then (i) the nanotube diameter
$d < 2 R \cos \theta$ and (ii) if the gap between vertically mesopores
is zero than the nanotube diameter must be less than
$d < \sqrt{4lR - l^2}$ (i.e.\ to avoid nanotube ``goiters'', see
further \cref{fig:upperbound_d}). The corresponding conditional
distribution is then given by
\begin{equation}
  \label{eq:Theta_d-given-R-theta-l}
  \Theta_d \given \Theta_R, \Theta_\theta, \Theta_l 
  \sim \unif(d_{-}, \min\{ \sqrt{4\Theta_l\Theta_R - \Theta_l^2},\, 
  2\Theta_R\cos\Theta_\theta,\, d_+ \})\,,
\end{equation}
where again we assume a uniform model. The form of the full
statistical model that includes the causal relationships in
\cref{fig:conditional-d-only} is then given by
\begin{equation}
  \label{eq:P_2}
  P_2 \defeq \delta_{U} (U - u(\B{x},t; \B{X})) \; 
  \delta_{\B{X}} (\B{X} -
  \B{\chi}(\B{\xi},t; \B{\Theta})) \; 
  \Pb(\Theta_d \given \Theta_R, \Theta_\theta, \Theta_l) 
  \Pb(\Theta_R) 
  \Pb(\Theta_\theta) \Pb(\Theta_l) \,,
\end{equation}
where the distributions for conditionals and priors above are given in
\cref{eq:Theta_R-Theta_t-Theta_l,eq:Theta_d-given-R-theta-l}.

We observe that the Bayesian networks in
\cref{fig:conditional-d-and-l} and \cref{fig:conditional-d-only} are
not equivalent. Over the physical range of hyperparameters in
\cref{tab:extended-range}, we observe that correlation structure in
\cref{fig:comparison-jkdes-corr-exAext} is distinct from
\cref{fig:comparison-jkdes-corr-exBext} and thus the joint density
$\Pb(\B{\Theta})$ for the Bayesian networks in
\cref{fig:conditional-d-only} and \cref{fig:conditional-d-and-l} are
different. Although it may be possible to select priors to ensure
equality between the distributions $P_1$ and $P_2$, this is not in
general the goal. Further, we can compare the Darcy-scale QoIs for the
model $P_2$ to the model $P_1$ in
\cref{fig:comparison-marginal-kdes-exAext} over the physical
hyperparameter ranges in \cref{tab:extended-range}, i.e.\ by comparing
\cref{fig:comparison-marginal-kdes-exAext} to
\cref{fig:comparison-marginal-kdes-exBext} and
\cref{fig:comparison-jkdes-corr-exAext} to
\cref{fig:comparison-jkdes-corr-exBext}. Importantly, incorporating
constraints directly into the models $P_1$ and $P_2$ allows us to
ensure that we sample from a joint distribution that ensures realistic
geometries that are consistent with the hierarchical nanoporous
material in \cref{fig:pore-structure} over the physical range of
hyperparameters. On the one hand, the densities in
\cref{fig:comparison-marginal-kdes-exBext} (and
\cref{fig:comparison-jkdes-corr-exBext}) are qualitatively similar to
those in \cref{fig:comparison-marginal-kdes-exAext} (respectively,
\cref{fig:comparison-jkdes-corr-exAext}) owing to the closeness of the
input densities $\Pb(\B{\Theta})$, see the empirical correlations in
\cref{fig:corr-exB-extended,fig:corr-exA-extended}. On the other hand,
we observe that the models $P_1$ and $P_2$ lead to distinct
distributions on Darcy-scale flow variables as the Cramér tests
reported in \cref{tab:cramer-test-univariate-P2} reject the hypothesis
of equality of the empirical distributions with high statistical
significance. We emphasize that we do not give here a method for
selecting among various models but instead a collection of tools for
breaking the stochastic modeling task into smaller, manageable
components that enables a systematic way of building a full
statistical model that incorporates engineering design constraints.

\begin{table}[!h]
  \centering
  \caption{Two-way nonparametric Cramér test
    (\cite{BaringhausFranz:2004tt}) rejects the hypothesis of equality
    of the empirical marginal and joint distributions for comparable
    variables in
    \cref{fig:comparison-marginal-kdes-exAext,fig:comparison-jkdes-corr-exAext,fig:comparison-marginal-kdes-exBext,fig:comparison-jkdes-corr-exBext}
    with high statistical significance (each test is based on
    \num{2000} values sampled using a gPCE).}
  \label{tab:cramer-test-univariate-P2}
  \begin{tabular}{lSSSS[table-format = <1.3]cl}
    \toprule
    Variable & \multicolumn{1}{c}{Cramér-statistic} & \multicolumn{1}{c}{Critical value} 
    & \multicolumn{1}{c}{Conf.\ interval} & \multicolumn{1}{c}{$p$-value} 
    & \multicolumn{1}{c}{Result} & Figures \quad \quad \\ % \midrule
    \midrule

    $\DL$& 4.688& 0.3389& 0.95& < 0.001 & reject
                                 & \rdelim\}{3}{0mm}[\ref{fig:comparison-marginal-kdes-exAext} vs.\ \ref{fig:comparison-marginal-kdes-exBext}]\\
    $\DT$& 79.08& 0.7396& 0.95& < 0.001 & reject & \\
    $\geff$& 0.4619& 0.3958& 0.95& 0.029 & reject & \\
    \midrule
    
    $(\DL,\DT)$& 72.31& 0.7827& 0.95& < 0.001& reject
                                 & \rdelim\}{3}{0mm}[\ref{fig:comparison-jkdes-corr-exAext} vs.\ \ref{fig:comparison-jkdes-corr-exBext}]\\
    $(\DT,\geff)$& 72.63& 0.8284& 0.95& < 0.001& reject & \\
    $(\geff,\DL)$& 4.734& 0.4882& 0.95& < 0.001& reject & \\ 
    \bottomrule
  \end{tabular}
\end{table}

\section{Conclusions}
\label{sec:conclusions}

The sensitivity of Darcy-scale observables to changes in pore-scale
properties and rigorous quantification of the uncertainty in
predictions are some of the least studied aspects of multiscale models
of flow and transport in porous materials. Our analysis leads to the
following major conclusions.
\begin{itemize}
\item Causal relationships are natural and stem from physical or
  chemical constraints, engineering design, and expert knowledge.
  These relationships exist between model parameters, scales, and
  model components in multi-scale and multi-physics models. In the
  context of hierarchical nanoporous materials, we observe that causal
  relationships suggested by geometrical structural constraints among
  microscopic features yield non-trivial correlations among pore-scale
  model parameters over physical ranges.
  
\item We incorporate causal relationships, and thereby correlations,
  in a unified random, multiscale PDE model using structured
  probabilistic graphical models, in this case Bayesian networks. This
  perspective ultimately gives rise to Bayesian network (random) PDEs.

\item Due to causal relationships and resulting correlations between
  model parameters, global sensitivity analysis is not
  straightforward. Furthermore, predictive PDFs of QoIs are not
  necessarily Gaussian or otherwise of known analytical form. For
  these two reasons, it is necessary to depart from moment-based
  sensitivity analysis, e.g. ANOVA methods, and deploy PDF-based
  methods such as mutual information.

\item The proposed mutual information global sensitivity indices for
  Bayesian network PDE and corresponding Darcy-scale QoIs yield
  parameter rankings that are consistent with our understanding of the
  physics of a hierarchical nanoporous material. We demonstrate that
  correlations stemming from causal relationships turn out to be
  important in determining the most influential model
  parameters/mechanisms and impact predictions of QoIs.

\end{itemize}

The structure of Bayesian networks or PGMs for the pore-scale
parameter correlations can in general be learned from experimental or
simulated data rather than based on natural structural constraints as
in this work and moreover they need not adopt causal relationships
(Bayesian networks) that mirror engineering design processes. One
future direction includes learning the structure of the Bayesian
networks from available data and inferring the distributions of
pore-scale features directly. The flexibility of the Bayesian network
PDE approach for incorporating uncertainty into physical models will
also enable explorations of model-form uncertainty that advance
traditional uncertainty quantification techniques. Another future
direction for research is to apply hybrid information divergences
(e.g.\ \cite{HallKatsoulakis:inprep}) to bound families of model
predictions based on various pore-scale parameter representations and
based on different upscaling techniques (i.e.\ different probabilistic
models) to address questions of model-form uncertainty and tracking
different levels of fidelity attached to system components.

\section*{Acknowledgments}

This research was supported in part by the Department of Energy (DOE)
Office of Science under grant 0000241868, by the Air Force Office of
Scientific Research under grants FA9550-17-1-0417 and
FA9550-18-1-0214, and by National Science Foundation under grant
1606192.

\bibliographystyle{siamplain}%
\bibliography{pgms_darcy_flow}%

\end{document}